\documentclass[12pt]{amsart}

\usepackage{amsfonts,amsmath,latexsym,amssymb,verbatim,amsbsy,times, esint}
\usepackage{amsthm}
\usepackage{pstricks}
\usepackage{mathtools}
\mathtoolsset{showonlyrefs}
\usepackage{graphicx}
\usepackage{caption}
\usepackage{tikz-cd}

\usepackage[top=0.8 in, bottom=0.8in, left=0.8in, right=0.8in]{geometry}

\usepackage[colorlinks=true, pdfstartview=FitV, linkcolor=blue,citecolor=green, urlcolor=blue]{hyperref}

\usepackage{enumitem}

\theoremstyle{plain}
\newtheorem{THEOREM}{Theorem}[section]

\newtheorem{LEMMA}[THEOREM]{Lemma}
\newtheorem{PROP}[THEOREM]{Proposition}

\theoremstyle{definition}

\theoremstyle{remark}
\newtheorem{REMARK}[THEOREM]{Remark}

\newcommand{\N}{\ensuremath{\mathbb{N}}}   
\newcommand{\R}{\ensuremath{\mathbb{R}}}   

\renewcommand{\S}{\ensuremath{\mathbb{S}}}

\def \d {\delta}

\def \e {\varepsilon}
\def \f {\varphi}

\def \k {\kappa}
\def \l {\lambda}
\def \L {\Lambda}

\def \n {\nabla}
\def \s {\sigma}
\def \th {\theta}

\def \w {\omega}

\def \id {{\bf id}}

\def \bn {{\bf n}}

\def \bu {{\bf u}}

\def\cB {\mathcal {B}}

\def \cD {\mathcal{D}}
\def \cE {\mathcal{E}}
\def\cF {\mathcal {F}}

\def \cH {\mathcal{H}}
\def \cJ {\mathcal{J}}
\def \cK {\mathcal{K}}
\def \cL {\mathcal{L}}

\def \cN {\mathcal{N}}

\def \cP {\mathcal{P}}
\def \cQ {\mathcal{Q}}
\def \cR {\mathcal{R}}

\def \p {\partial}

\DeclareMathOperator{\Span}{span} %
\DeclareMathOperator{\supp}{supp} %
\DeclareMathOperator{\diver}{div} %

\def \cP {\mathcal{P}}

\def \dd  {\mathrm{d}}

\def \ds  {\, \mathrm{d}s}
\def \dt  {\, \mathrm{d}t}
\def \dtau  {\, \mathrm{d}\tau}

\def \dx  {\, \mathrm{d}x}

\def \dr  {\, \mathrm{d}r}

\def \dr  {\, \mathrm{d}r}
\def \dR  {\, \mathrm{d}R}
\def \ds  {\, \mathrm{d}s}

\def \sfS  { \mathsf{S}}
\def \sfT  { \mathsf{T}}

\def \brho {\bar{\rho}}

\def \temp {\mathrm{temp}}
\def \tot {\mathrm{tot}}

\def \eff {\mathrm{eff}}

\usetikzlibrary{cd}
\usetikzlibrary{snakes}

	\title[Stability for Hyperelastic Bodies]{Stability of Self-Similar Expansions \\ for Self-Gravitating Hyperelastic Bodies}
	
	\author{Juhi Jang and Trevor M. Leslie}
	
\begin{document}

\vspace*{-5 mm}
\begin{abstract}
We consider the three-dimensional elastodynamics system with nonlocal gravitational coupling and nonlinear traction-free boundary conditions, modeling an isolated self-gravitating hyperelastic body, where the strain-energy function is assumed to be objective, isotropic, and homogeneous of degree~$-1$. We prove nonlinear stability of a self-similarly expanding separable solution, recently constructed by the authors \cite{JangLeslie2025}, under radial perturbations satisfying the codimension-one constraint of zero total energy. Using a self-similar reformulation analogous to that developed for the mass-critical Euler--Poisson system in \cite{HadzicJang2018CPAM}, we recast the perturbed equation as a scalar nonlinear five-dimensional radial wave equation and the nonlinear traction-free boundary condition as a homogeneous Neumann condition. A key is to identify a quasilinear form of the leading-order spatial operator that is compatible with this boundary condition and yields coercive energy estimates. Our energy estimates capture the stabilizing effect of the background expansion and yield global existence and nonlinear codimension-one stability in self-similar variables. 
\end{abstract}

\maketitle

\begin{small} 
\tableofcontents
\end{small} 

\vspace*{-5 mm}
\section{Introduction}

Consider the following equations of motion for a self-gravitating elastic body $\cB(t)$ in vacuum.
\begin{equation}
	\label{e:Eulerian}
	\begin{cases} 
		\p_t \rho(t,x) + \diver(\rho \bu)(t,x) = 0,
		\qquad x\in \cB(t) = \supp \rho(t), \\
		\rho (\p_t\bu + \bu\cdot \n \bu)(t,x) = \diver \sfT(t,x) - G\rho(t,x) \n_x \Phi[\rho(t)](x), \qquad x\in \cB(t),\\
		\Delta \Phi[\varrho] = 4\pi \varrho, \qquad \lim_{|x|\to \infty} \Phi[\varrho](x) = 0 \\
		\sfT(t,x)\bn(x) = 0, \qquad x\in \p \cB(t).
	\end{cases}
\end{equation}
The functions $\rho\ge 0$ and $\bu$ correspond to the density and velocity of the body, respectively, while $\bn(x)$ is the unit normal vector to the boundary of $\cB(t) = \supp \rho(t)$, and $G$ is the gravitational constant.  Finally, $\sfT$ denotes the \textit{Cauchy stress tensor}; one must specify how $\sfT$ responds to the dynamics (a \textit{constitutive relation}) in order to close the system \eqref{e:Eulerian}.

This paper is concerned with the stability of solutions to \eqref{e:Eulerian}.  More specifically, we are concerned with the codimension-one stability of a particular class of radially symmetric, separable, self-similar expanding solutions (which we refer to as the `background solution') that we specify later on.  Since $\cB(t)$ evolves in time, a naive comparison between Eulerian configurations is impractical, and our statement and proof are formulated in Lagrangian coordinates.   At a high level, our analysis is similar in approach to that of the first author's joint work \cite{HadzicJang2018CPAM} with Hadžić on the mass-critical Euler--Poisson equations. However, our new stress tensor $\sfT$ and boundary conditions $\sfT \mathbf{n} = 0$ (the traction-free solid boundary conditions), introduce significant structural changes in the system that have prompted us to substantially rebuild the technical framework of Hadžić and Jang's original argument.  Furthermore, part of our contribution in this work is to elucidate some of the key structures of the equation for the perturbation around the background solution.  We argue that, under our constitutive assumptions, the equation for the perturbation takes the form of a nonlinear wave equation 
with pure Neumann boundary conditions.  The simple form of our boundary conditions is especially remarkable, given that the boundary conditions for the full system are nonlinear and nonlocal, even under symmetry assumptions.

The rest of the Introduction is organized as follows.  First, we set the stage with some background information: We establish notation and terminology; we specify our constitutive assumptions; and we  provide the explicit equations (in Lagrangian coordinates) for radial and separable dynamics under our assumptions.  Many of our conventions are similar to those of our earlier paper \cite{JangLeslie2025}, and when discussing such issues we refer to the latter for more details and motivation, keeping the presentation here self-contained but as streamlined as possible. Next, we explicitly record the system of equations satisfied by the perturbation itself, and we define a hierarchy of perturbation energies that provide us with an appropriate framework to study the stability of our perturbations.  Finally, we state our main results and contextualize them within the existing literature. 

\subsection{Dynamics and constitutive assumptions}

\subsubsection{Lagrangian coordinates and spherical symmetry}
As mentioned above, we work in Lagrangian coordinates. Our reference domain will be $\cB = \{X\in \R^3: |X|\le 1\}$. Letting $|\cB|$ denote the Lebesgue measure of $\cB$, we define our (constant) reference density $\brho$ via $\brho = \frac{1}{|\cB|}\int_{\cB(0)} \rho(0,x) \dx$.
We use $x:[0,T)\times \cB \to \cB(t)$ to denote the flow map associated to $(\rho, \bu)$, satisfying  $\dot{x}(t,X) = \bu(t,x(t,X))$ and $\rho(t, x(t,X)) \det Dx(t,X) = \brho$, for all $t$ belonging to some maximal time interval $[0,T)$ and all $X\in \cB$.  As a baseline assumption, the flow map will always be at least $C^2$ in time, and $x(t, \cdot)$ will be an orientation-preserving $C^2$ diffeomorphism for each $t\in [0,T)$.  We consider only spherically symmetric solutions, $x(t, X) = \phi(t, |X|) \frac{X}{|X|}$ for $X\ne 0$ and some function $\phi:[0,T)\times [0,1]\to [0,+\infty)$, so the assumption that $x(t, \cdot)$ is an orientation-preserving $C^2$ diffeomorphism forces $\phi(t,0) = \phi_{RR}(t,0)~=~0$ and $\phi_R(t,R)>0$ (cf. \cite{Sideris2023}, Section 5).  In what follows, we will refer to $\phi$ as a `solution', with the understanding that $\rho$ and $\bu$ can be recovered from $\phi$ and its derivatives.  The eigenvalues of the deformation tensor $Dx(t,X)$ are $\phi_R(t,R)$, which has multiplicity one, and 
\begin{equation} 
	\label{e:Lambdadef}
	\L(t,R) = \begin{cases} 
		\frac{\phi(t,R)}{R} & R\in (0,1], \\
		\phi_R(t,0) & R = 0,
	\end{cases}
\end{equation}
which has multiplicity two. The ratio of these eigenvalues, 
\begin{equation}
	\label{e:ydef}
	y(t,R) = \frac{\phi_R(t,R)}{\L(t,R)}, 
\end{equation}
will appear frequently in our analysis.

\subsubsection{Separable solutions and perturbations}

The term \textit{separable solution} will refer to a spherically symmetric  solution (which we will usually distinguish from other spherically symmetric solutions by writing $\phi_*$ instead of $\phi$) that takes the form $\phi_*(t,R) = q(t)\f(R)$, where $q:[0,T)\to [0,\infty)$ and $\f:[0,1]\to [0,\infty)$ are both $C^2$.  We will always normalize $q$ so that $q(0) = 1$. Separable solutions are sometimes referred to as \textit{homologous} by other authors.  Under appropriate constitutive assumptions, the authors have constructed a family of separable solutions in \cite{JangLeslie2025}, following the strategy of Sideris's earlier paper \cite{Sideris2023}, which treated the case $G = 0$.  See also \cite{Calogero2022} for another construction using different constitutive assumptions and different techniques.  

This paper is concerned with understanding the stability of spherically symmetric perturbations of separable solutions like those constructed in \cite{JangLeslie2025}.  More specifically, we will deal with zero-energy perturbations of self-similar expanding separable solutions; we will explain below what we mean by `self-similar' in this context.  We make the ansatz
\begin{equation}
	\label{e:ansatz}
	\phi(t,R) = \phi_*(t,R) (1 + \xi(t,R)), 
	\qquad \phi_*(t,R) = q(t) \f(R),
\end{equation}
and set notation analogous to \eqref{e:Lambdadef}, \eqref{e:ydef}:
\begin{equation}
	\label{e:deflambday*z}
	\l(R) = \begin{cases}
		\frac{\f(R)}{R}, & R\in (0,1], \\
		\f'(0), & R = 0, 
	\end{cases}
	\qquad 
	y_*(R) = \frac{\f'(R)}{\l(R)}, 
	\qquad 
	z(t,R) = \frac{R \xi_R(t,R)}{1 + \xi(t,R)}.
\end{equation}
In this notation, we have $y(t,R) = y_*(R) + z(t,R)$.  We will refer to $\phi_*(t,R)$ as the \textit{background solution}, and we will call $\phi(t,R)$ a \textit{perturbed} solution with \textit{perturbation} $\xi(t,R)$.

\subsubsection{Piola--Kirchhoff stress and constitutive assumptions}
Given a flow map $x(t,X)$ with deformation tensor $Dx(t,X)$ and Jacobian determinant $Jx(t,X) = \det Dx(t,X)$, the \textit{Piola--Kirchhoff stress tensor} $\sfS$ associated to  $(\cB, \brho)$ is defined via
\begin{equation}
	\label{e:Sdef}
	\sfS(t,X) = Jx(t,X) \sfT(t, x(t,X)) Dx(t,X)^{-T},
	\qquad (t,X)\in [0,T)\times \cB.
\end{equation}
Our constitutive assumptions, listed presently, will be stated in terms of $\sfS$ and related functions.  In our list of assumptions, $GL_+(3, \R)$ denotes the space of invertible $3\times 3$ matrices with real entries that have positive determinant, and $SO(3,\R)$ is its subgroup consisting of $3\times 3$ orthogonal matrices with determinant $1$.  Most of our constitutive assumptions correspond to well-established terms in the continuum mechanics literature; we point out those terms where appropriate.
\begin{enumerate}[label=(\alph*)]
	\item We assume $\sfS$ takes the form
	\[
	\sfS(t,X) = \widehat{\sfS}(Dx(t,X),X), 
	\quad 
	\widehat{\sfS}(F,X) = \brho \frac{\p}{\p F}W(F,X),
	\quad F\in GL_+(3,\R), \;X\in \cB, \; t\in [0,T),
	\]
	where $\widehat{\sfS}$ is at least $C^1$ and $W$ is at least $C^2$. The existence of a `constitutive function' $\widehat{\sfS}$ says that the body under consideration is \textit{elastic}; the existence of a \textit{strain-energy} function $W$ says that the body is \textit{hyperelastic}.
	\item We assume that 
	\[
	W(OF,X) = W(F,X) = W(FO, X), \qquad F\in GL_+(3,\R), \;O\in SO(3, \R).
	\]
	The first equality says that $W$ is \textit{objective}; the second says that it is \textit{isotropic}.
	
	\item We assume that $W(F,X) = W(F)$ is independent of $X$.  (This says that the body is \textit{homogeneous}.)
	
	\item We assume that the function $W$ possesses homogeneity of degree $-1$:
	\[
	W(\s F) = \s^{-1}W(F),
	\qquad \s>0, \;F\in GL_+(3,\R).
	\]
	For definiteness, we also fix the normalization $W(I) = \brho^{\frac13}$.
\end{enumerate}
Assumptions (a)--(c) are rather mild, while (d) is a strong hypothesis. Some homogeneity condition, however, is needed in order to produce separable solutions, cf. \cite{Sideris2023, JangLeslie2025}, and in the presence of gravity (or more specifically, the scaling of the Newtonian potential term), only degree-$(-1)$ homogeneity is compatible with separable solutions.  Since we are concerned with the stability of such solutions, assumption (d) is necessary for our framework.  

Assumptions (a)--(d) guarantee that the association $(0,\infty) \times \S^2 \ni (y,\w)\mapsto W\big(I + (y-1)(\w \otimes \w)\big)$ is actually independent of $\w$.  (Here $\S^2$ is the unit $2$-sphere in three dimensions.) Therefore, we can define a function $g:(0,\infty)\to \R$ via 
\begin{equation} 
	\label{e:gdef} 
	g(y) = \brho^{- \frac13} W\big(I + (y-1)(\w \otimes \w)\big),
	\qquad \w\in \S^2.
\end{equation} 
By design (and because of the normalization of $W$), we have $g(1) = 1$, and our constitutive assumptions then force $g'(1) = -\frac13$.  We generally work in the setting where $|y-1|$ is small and assume that $g$ is sufficiently regular for our computations to go through, with $g''(1)>0$ (and therefore $g''(y)>0$ if $|y-1|$ is small enough).  Note that the choice $g(y) = y^{-\frac13}$ gives rise to the Cauchy stress tensor of a mass-critical polytropic gas (i.e., a gas with associated pressure $p = \frac13\rho^{\frac43}$).  Our `solid boundary' conditions \eqref{e:Eulerian}$_4$ are different from the ones typically used for a polytropic gas, and our theory does not apply directly to this case, but this example is nonetheless useful to keep in mind as a quick plausibility check on the numerology.
In fact, our constitutive assumptions (and in particular the degree-$(-1)$ homogeneity of $W$) also give rise to the same scaling symmetry that is exhibited by a mass-critical polytropic gas: If $x(t,X)$ is a solution of the Lagrangian formulation of \eqref{e:Eulerian}, then so is 
\begin{equation}
\label{e:scaling}
x_a(t,X) = a x(a^{-\frac32}t, X),
\end{equation} 
with the same reference configuration $\cB$ and reference density $\brho$.  If $(\rho, \bu)$ is the density/velocity pair associated to $x(t,X)$, then this means that $(\rho_a, \bu_a)$ (together with an appropriately rescaled deformation tensor) also satisfies \eqref{e:Eulerian}, where  
\begin{equation}
\rho_a(t,x) = a^{-3} \rho(a^{-\frac32} t, a^{-1}x), 
\qquad 
\bu_a(t,x) = a^{-\frac12} \bu(a^{-\frac32}t, a^{-1}x), 
\qquad \qquad a>0.
\end{equation} 
These scaling symmetries are intimately related to the existence of the separable solutions, the temporal and spatial dynamics of which are completely decoupled, and the stability of which we study in this paper.

\subsubsection{Equations for spherically symmetric and separable solutions} It is shown in \cite{JangLeslie2025} that, under the above assumptions, the system \eqref{e:Eulerian} becomes the following, for spherically symmetric solutions. Here, $y = y(t,R)$ is as in \eqref{e:ydef}.
\begin{equation}	
	\label{e:spherical}
	\begin{cases}
		\displaystyle \phi_{tt} - \brho^{\frac13} \bigg[\frac{1}{R^2\phi} \p_R \bigg( \frac{R^4}{\phi} {g}'({y}) \bigg) + \frac{R}{\phi^2} {g}({y})\bigg] +  \frac{4\pi}{3} R^3 G \brho \phi^{-2} = 0, & \qquad \text{ in } [0,T)\times \cB, \\
		{g}'({y}(t,1)) = 0, \\
		\phi(t,0) = \phi_{RR}(t,0) = 0,\; \phi_R(t,\cdot)>0.
	\end{cases}
\end{equation} 

For \textit{separable} solutions $\phi_*(t,R) = q(t) \f(R)$ (with $q(0) = 1$ by assumption), the above simplifies even more.  Denote $y_* = y_*(R)$ as in \eqref{e:deflambday*z}.  There exists $\mu\in \R$ such that 
\begin{equation}
	\label{e:separated}
	\begin{cases}
		\displaystyle 
		\mu \f - \brho^{\frac13} \bigg[ \frac{1}{R^2 \f} \p_R \bigg( \frac{R^4}{\f} g'(y_*)\bigg) + \frac{R}{\f^2}g(y_*)\bigg] + \frac{4\pi}{3}R^3 G \brho \f^{-2} = 0, \\
		g'(y_*(1)) = 0, \;\f(1) = 1, \\
		\f(0) = \f''(0) = 0, \; \f' >0.
	\end{cases}
\end{equation}
Furthermore, 
\begin{equation}
	\label{e:q}
	q^2 \ddot{q} = \mu,
	\;\; q(0) = 1.
\end{equation}	

Given a background solution $\phi_*(t,R) = q(t) \f(R)$ (with $\varphi$ and $q$ satisfying \eqref{e:separated} and \eqref{e:q}, respectively) and a perturbed solution $\phi(t,R) = \phi_*(t,R)(1+\xi(t,R))$ satisfying \eqref{e:spherical}, we can derive an equation for the  perturbation $\xi$; see Section \ref{ss:Equations} below.  Note that our definition of $\l$, together with the requirement that $\f'>0$,  ensures that $\l$ is bounded above and below by positive constants.

\subsubsection{Expansion and collapse of separable solutions: the $q$ equation}
Equation \eqref{e:q} is elementary \cite{FuLin1998, HadzicJang2018CPAM}, but its solutions have some very useful and relevant properties that we pause to record.  First of all, a function $q(t)$ satisfying \eqref{e:q} has an associated conserved quantity, which we call the \textit{effective energy}:
\[
\frac12 \dot{q}(t)^2 + \frac{\mu}{q(t)} = \frac12 \dot{q}(0)^2 + \mu =:e_{\eff}.
\]
The signs of $\mu$, $e_{\eff}$, and $\dot{q}(0)$ determine the eventual fate of $q(t)$ and therefore of a separable solution $\phi_*(t,R) = q(t)\f(R)$.  If $\mu>0$, then $q(t)$ eventually grows at an asymptotically linear rate (the `linearly expanding' case).  If $\mu=0$, then $\dot{q}(t)\equiv \dot{q}(0)$, so we have either linear expansion, linear collapse, or a steady state.  The most interesting situation is the one where $\mu<0$, and it is compatible with several qualitatively distinct behaviors.  Assume, therefore, that $\mu<0$.  If $e_{\eff}$ and $\dot{q}(0)$ are both positive, then $q(t)$ exists for all time and experiences asymptotically linear expansion as $t\to \infty$.  If $e_\eff = 0$, then $q(t)$ is given by the formula $q(t) = (1 + \frac32 \dot{q}(0)t)^{\frac23}$ for all $t$ in its domain of existence, which is finite if $\dot{q}(0)<0$ and includes all of $[0,\infty)$ if $\dot{q}(0)>0$.  We refer to this situation as either `self-similar expansion' or `self-similar collapse', depending on the sign of $\dot{q}(0)$, on the grounds that in this case, the natural rescaling \eqref{e:scaling} results in $q_a(t):=a q(a^{-\frac32} t) = q(t - t_a)$, with $t_a = \frac23\dot{q}(0)^{-1}(1-a^{\frac32})$; that is, the rescaled $q_a(t)$ is a translation of the original function $q(t)$, for these solutions, and the same is true of $x_a$ and $x$ from \eqref{e:scaling}.  All other configurations of $e_\eff$ and $\dot{q}(0)$ result in finite-time collapse, with the same `self-similar' rate, i.e., $q(t){\sim}_{t\to T^-} \;(T - t)^{\frac23}$, where $T$ is the time of collapse.  

Our focus will be on the case of self-similar expansion.  Within the class of separable solutions satisfying $\mu<0$, $\dot{q}(0)>0$, this situation lies on the border $e_{\eff}=0$ between finite-time collapse ($e_\eff<0$) and asymptotically linear expansion ($e_{\eff}>0$).  So it is already clear that self-similar expanding solutions are not stable outright, even within the class of separable solutions.  We claim, however, that a codimension-one stability holds for these solutions.  We consider \textit{zero-energy} perturbations of self-similar expanding solutions.  We spell out what this means below.

\subsubsection{Energy}

Under our constitutive assumptions (hyperelasticity being the most salient), the system \eqref{e:Eulerian} has the following natural energy.  
\begin{equation}
	\label{e:Energy1}
	E(\rho, \bu)(t)
	= \frac12 \int_{\cB(t)} \rho(t,x)\left( |\bu(t,x)|^2 + G\Phi[\rho(t)](x) \right) \dx + \int_{\cB} \brho W(Dx(t,X)) \,\dd X. 
\end{equation} 
A quick computation shows that energy is conserved: $\frac{\dd}{\dt} E(\rho, \bu)(t) = 0$.  For spherically symmetric solutions, it is also possible to rewrite this energy in the form
\begin{equation} 
	\label{e:Energyspherical}
	E(\rho, \bu) = E(\phi) = \int_0^1 4\pi R^2 \brho \bigg[ \frac12 \phi_t^2 + \phi \phi_{tt} \bigg] \,\dd R.
\end{equation} 
It would of course be more explicit to replace  $\phi_{tt}$ in this expression using \eqref{e:spherical}$_1$, but in the present form, it is easy to see what happens when $\phi(t,R) = \phi_*(t,R) = q(t)\f(R)$ is a separable solution, in which case the quantity in brackets is simply $e_{\eff}\cdot \f^2$. The energy of a separable solution is thus given by 
\begin{equation}
\label{e:Energyseparable}
	E(\rho_*, \bu_*) = E(\phi_*) = e_{\eff} \int_0^1 4\pi \brho \l(R)^2 R^4 \, \dd R.
\end{equation}
In particular, a separable solution has zero energy if and only if its effective energy $e_{\eff}$ is zero.  Given a separable solution $\phi_*(t,R) = q(t)\f(R)$ with zero energy, what we mean by a `zero-energy perturbation' is a perturbation $\xi(t,R)$ such that the corresponding perturbed solution $\phi(t,R) = \phi_*(t,R)(1 + \xi(t,R))$ still lies on the codimension-one manifold of solutions such that $E(\rho, \bu) = E(\phi) = 0$.  Note that below we will also define a different kind of `energy', not to be confused with $E(\rho, \bu)$, that we will refer to as \textit{perturbation energy}.  The latter will provide us with a tool for making precise the sense in which our solutions are stable. 

\subsection{Equations and energy hierarchies for the perturbation}

\label{ss:Equations}

Unless otherwise stated, we assume from now on that $\phi$ is a solution to \eqref{e:spherical} satisfying the ansatz \eqref{e:ansatz}, with $\phi_*(t,R) = q(t)\f(R)$ a separable background solution associated to $\phi$ and $\xi(t,R)$ the associated perturbation.  In this subsection, we record the equations relevant to the dynamics and energy of the perturbation $\xi$.  

\subsubsection{Equations of motion}

For a solution $\phi$ satisfying the ansatz \eqref{e:ansatz}, we obtain an equation for the perturbation $\xi$ as follows.  We substitute $\phi(t,R) = q(t)\f(R)(1 + \xi(t,R))$ into \eqref{e:spherical}$_1$ and multiply the result by $\frac{q^2}{\f}$; then we add $-\frac{1}{\f(1+\xi)^2}$ times \eqref{e:separated}$_1$.  The gravitational terms cancel, and after a few manipulations, we obtain 
\begin{equation}
	\label{e:xit}
		q^2 (q \xi)_{tt} + 2 \mu \xi + \cL_{\xi}(\xi) = \cN, 
\end{equation}
where 
\begin{equation} 
\label{e:Lmuxi}
	\cL_{\xi}(u) = -\frac{\brho^{\frac13}}{R^4 \l^2 (1+\xi)} \p_R \bigg(\frac{1}{\l (1+\xi)^2} \cdot  \frac{g'(y_*+z) - g'(y_*)}{z} \cdot   R^4 u_R \bigg),
\end{equation} 
\begin{equation}
	\cN = \xi^2 \cN_1 + \xi_R^2 \,\cN_2,
\end{equation}
\begin{equation}
	\cN_1 = - \frac{\mu}{\xi^2} \bigg( 1 - 2\xi - \frac{1}{(1+\xi)^2} \bigg) = \bigg[ 2 + \frac{1}{1+\xi} \bigg] \frac{\mu}{1+\xi},
\end{equation} 
\begin{equation}
\label{e:N2}
	\cN_2 = \frac{\brho^{\frac13}}{\l^3 (1+\xi)^4} \cdot \frac{g(y_*+z) - g(y_*) - z g'(y_*)}{z^2}. 
\end{equation}

As for the boundary conditions, note that \eqref{e:separated}$_2$ requires $g'(y_*(1)) = 0$, while \eqref{e:spherical}$_2$ together with our ansatz requires that $g'(y(1)) = g'(y_*(1) + z(t,1)) = 0$.  In our regime where $\xi$, $\xi_R$, and $y_*-1$ are small and $g''>0$ near $1$, these two conditions are equivalent to requiring that $z(t,1) = 0$, which in turn is equivalent to $\xi_{R}(t,1) = 0$.  Similarly, combining \eqref{e:spherical}$_3$ and \eqref{e:separated}$_3$ yields $\xi_R(t,0) = 0$, after a short calculation.  In summary, the boundary conditions for $\xi$ reduce to simple \textit{Neumann} boundary conditions rather than the nonlinear boundary conditions of \eqref{e:spherical}:
\begin{equation} 
\label{e:xiBC}
\xi_R(t,0) = \xi_R(t,1) = 0,
\qquad t\ge 0.
\end{equation} 

\begin{REMARK}
\label{r:Lform}
The quantity $\cN$ consists of terms that are essentially quadratic, while $\cL_{\xi}(\xi)$ contains all the linear contributions from the spatial derivatives.  However, $\cL_{\xi}(\xi)$ also contains many nonlinear terms, some of which can easily be moved to the right hand side of \eqref{e:xit} if desired.  The reason we settle on this particular operator is because its desirable prefactor $\frac{1}{R^4 \l^2 (1+\xi)}$ guarantees that it contributes nothing to the energy $E(\phi)$.  Indeed, if we make the ansatz $\phi(t,R) = q(t)\f(R)(1 + \xi(t,R))$, substitute into \eqref{e:Energyspherical}, and use \eqref{e:xit}, it results in the appearance of the term $R^4 \l^2 (1+\xi)\cL_{\xi}(\xi)$ (and no other terms involving $\cL_{\xi}(\xi)$).  This term vanishes upon integration, by design.  
\end{REMARK}

\subsubsection{Energy for the perturbed solution} Let us record the energy of a perturbed solution in terms of the energy of a separable background solution.  

We substitute the ansatz $\phi(t,R) = q(t) \f(R) (1+\xi(t,R))$ into \eqref{e:Energyspherical}. Using $q^2 q_{tt} = \mu$ then yields
\begin{equation}
	E(\rho, \bu)(t) = 4\pi \brho \int_0^1 \bigg[ \frac{1}{2} \big(q(1+\xi)\big)_t^2 + \frac{1+\xi}{q}\big( \mu + q^2 (q\xi)_{tt} \big)  \bigg] \l^2 R^4 \,\dd R. 
\end{equation}
We use \eqref{e:xit} to replace $q^2(q\xi)_{tt}$ in the integral above.  As remarked above, the operator $\cL_{\xi}(\xi)$ is designed precisely so that its contribution to the energy will vanish here.  Next, use $\frac12 q_t^2 + \frac{\mu}{q} = e_{\mathrm{eff}}$ and \eqref{e:Energyseparable}, to conclude that 
\begin{equation} 
	\label{e:EnergyPerturbed}
	E(\phi)
	= E(\phi_*) + 4\pi \brho \int_0^1 \bigg[ \bigg( 2e_{\eff} - \frac{3\mu}{q}\bigg) \xi + q_t q \xi_t \bigg] R^4 \l^2 \dR - \frac{4\pi \brho}{q}\cJ(t), 
\end{equation}
where
\begin{equation}
	\label{e:Jt}
	\cJ(t) = - \int_0^1 \bigg( \frac12 q(q\xi)_t^2 + \frac{\mu \xi^2}{1 + \xi} + (1 + \xi)\xi_R^2\,\cN_2 \bigg) R^4 \l^2 \dR
\end{equation}
encodes the terms that are nonlinear in $\xi$ and its derivatives.  

\subsubsection{Zero-energy self-similar expanding solutions}

As mentioned above, our focus in this paper is on the codimension-one stability of self-similarly expanding solutions.  This has two consequences for our setup.  First, the case we are interested in is the one where $\mu<0$, $e_{\eff} = 0$, and $\dot{q}(0)>0$, in which case $q(t) = (1 + \frac32 \dot{q}(0)t)^{\frac23}$, and it is appropriate to rescale time through the new variable 
\begin{equation}
\label{e:timescales}
s(t) = \int_0^t \frac{\dt'}{q(t')^{\frac32}} = \frac{1}{\frac32\dot{q}(0)} \log\big( 1 + \frac32 \dot{q}(0) t \big),
\qquad t\ge 0.
\end{equation} 
This causes $\frac{q_s}{q} = \sqrt{-2\mu}$ to be constant; we set the following notation: 
\[
b  :=  -\sqrt{-2\mu}<0
\quad\qquad \text{(so $\mu = -\tfrac12 b^2$).}
\]
In the new variables, the two equations \eqref{e:xit}, \eqref{e:xiBC} become the system 
\begin{equation}
	\label{e:xis}
	\begin{cases}
		\displaystyle
		\xi_{ss} - \frac12 b\xi_s + 3\mu \xi + \cL_{\xi}(\xi) = \cN, \\[1ex]
		\xi_{R}(s,0) = \xi_{R}(s,1) = 0.
	\end{cases}
\end{equation}

The other consequence of our codimension-one stability setup is that the energy of the perturbative solution we are interested in must also be zero.  Setting $E(\phi)$ and $E(\phi_*)$ both equal to zero in \eqref{e:EnergyPerturbed}, and writing the result (as well as \eqref{e:Jt}) in terms of the new time variable $s$, we obtain 
\begin{equation}
	\label{e:Energys}
	\int_0^1 \bigg( \frac32 b^2 \xi - b \xi_s \bigg) R^4 \l^2 \,\dR = \cJ(s),
\end{equation}
\begin{equation}
	\label{e:Js}
	\cJ(s) = - \int_0^1 \bigg( \frac12 |\xi_s - b\xi|^2 + \frac{\mu \xi^2}{1 + \xi} + (1 + \xi)\xi_R^2\,\cN_2 \bigg) R^4 \l^2 \dR.
\end{equation}

\subsubsection{Perturbation energies}

Let us introduce the operator
\begin{equation} 
	\label{e:DRdef}
	D_R = \frac{1}{R^4} \p_R \big( R^4 \cdot \big) = \p_R + \frac4R,
\end{equation} 
which we will use repeatedly below.  We also define $\cD^k$ via $\cD^0 = \id$, $\cD^1 = \p_R$, $\cD^2 = D_R \p_R$, $\cD^3 = \p_R D_R \p_R$, etc.  We define presently a hierarchy of perturbation energies for the system \eqref{e:xis} in terms of these operators $\cD^k$.  We will outline some properties of the calculus of these operators in Section \ref{ss:Preliminaries}.

\label{sss:PerturbationEnergies}

We define the `total' perturbation energy at levels zero through $j$ via 
\begin{equation} 
	\label{e:Etotdef}
	\cE_{\tot}^{\le j}(s) = \sum_{\substack{\ell+m\le j \\ m \ge 0;\;\;\ell \ge 0}} \int_0^1 |\p_s^\ell \cD^m \xi(s,R)|^2 R^4 \dR,
\end{equation}
and we define the corresponding `temporal' perturbation energy similarly, but restricting the number of spatial derivatives to zero or one:
\begin{equation}
	\label{e:Etempdef}
	\cE_{\temp}^{\le j}(s) = \sum_{\substack{\ell+m\le j \\ m \in \{0,1\};\;\;\ell \ge 0}} \int_0^1 |\p_s^\ell \p_R^m \xi(s,R)|^2 R^4 \dR.
\end{equation}

Our stability result, loosely stated, says that if the perturbation energy $\cE^{\le J}_{\tot}$ of some valid initial condition is small enough initially (for some $J\ge 5$), then this initial condition will give rise to a global-in-time solution, for which $\cE^{\le J}_{\tot}$ shrinks to zero asymptotically.  (We will specify below what we mean by `valid' initial condition.) Though the temporal energy does not appear explicitly in this statement, it plays a key role in the proof.  It is clear that $\cE^{\le j}_{\temp} \le \cE^{\le j}_{\tot}$ for any $j$.  However, if $\xi$ solves \eqref{e:xis}, then, as we show in Section \ref{s:EquivEnergies}, it turns out that $\cE^{\le J}_{\tot}\lesssim \cE^{\le J}_{\temp}$ when $J \ge 5$ and when $\cE_{\tot}^{\le J}$ is small enough.  This will allow us to interchange $\cE_{\temp}^{\le J}$ with $\cE^{\le J}_{\tot}$ in a key place in the proof of our main stability estimate.  We formally record this `equivalence of energies' statement in Lemma \ref{l:EnEquiv} below.

\subsubsection{Initial conditions and compatibility conditions} 
\label{sss:IC}
As we have just mentioned, our existence and stability result relies on a smallness condition on $\cE^{\le J}_{\tot}$ at time zero.  We will use $\cE_0^{\le J}$ to denote the initial total energy up to level $J$, associated to a set of initial conditions $(\xi_0, \xi_1)$ for $(\xi, \xi_s)$.  However, we need access to the initial conditions of higher-order time derivatives in order to define this properly.  Note that the function $\partial_s^j \xi(0,\cdot)$ is determined from the lower-order time derivatives $(\partial_s^i \xi(0,\cdot))_{i=0}^{j-1}$ through the equation \eqref{e:xis}$_1$, for $j=2, \ldots, J$ (and the background solution). We may therefore let $(\xi_j)_{j=2}^J$ denote the initial conditions associated to $(\partial_s^j \xi)_{j=2}^J$; they are determined completely from $\xi_0$ and $\xi_1$.  Note that if the boundary conditions $\xi_R(s,1)=0$ are satisfied, then $\xi_j$ will need to satisfy the compatibility  $\partial_R \xi_j(1) = 0$ for each $j\in \{0, \ldots, J-2\}$.  This is really a restriction on $\xi_0$ and $\xi_1$; we will therefore say that the initial condition $(\xi_0, \xi_1)$ is `compatible with the boundary conditions of the equation' if $\xi_0, \ldots, \xi_{J-2}$ satisfy the compatibility conditions.  
With the above notation in hand, we may define $\cE_0^{\le J}$ via
\begin{equation}
\label{e:initialtotalenergy}
\cE_0^{\le J}  = \sum_{\substack{\ell + m \le J \\ m\ge 0, \; \ell \ge 0}} \int_0^1 |\cD^m \xi_\ell(R)|^2 R^4\dR.
\end{equation} 
Finally, we stress once again that our stability theory for self-similar solutions is predicated on the assumption that the perturbations we consider satisfy the zero-energy conditions \eqref{e:Energys}--\eqref{e:Js}.  Therefore, for purposes of our stability estimates, we need \eqref{e:Energys}--\eqref{e:Js} to be satisfied with $(\xi_0, \xi_1)$ in place of $(\xi(s), \xi_s(s))$.  

\subsection{Main result, methodology, and related literature}

We are now in a position to state the main result of this paper, which establishes   global nonlinear codimension-one stability of self-similarly expanding finite self-gravitating elastic bodies under zero-energy radial perturbations. To the best of our knowledge, this theorem provides the first  global-in-time solutions for self-gravitating hyperelastic balls beyond the separable class. 

\begin{THEOREM}
\label{t:main}
Assume $g\in C^{J+2}((0,\infty))$, $J\ge 5$, and that $g(1) = 1$, $g'(1) = -\frac13$, and $g''(1)>0$.  There exists $\e>0$, depending only on $g$ and $J$, and an absolute constant $c>0$, such that the following holds.  

Given $\brho>0$ and $\mu\in (-c\brho^{\frac13}g''(1),0)$, let $\varphi$ be a classical solution of \eqref{e:separated} (equipped with the parameters $\brho$, $\mu$) that satisfies  $|\varphi''(R)|\le R\e$, $R\in (0,1)$.  Then there exist positive constants $\theta$, $\kappa$, and $C$, depending on $g$, $J$, $\varphi$, $\mu$, and $\brho$, with the following property.  If $(\xi_0, \xi_1)$ is a set of initial conditions that is compatible with the boundary conditions of \eqref{e:xis} and satisfies the zero-energy condition \mbox{\eqref{e:Energys}--\eqref{e:Js}}, and if  $\cE_0^{\le J} \le \theta$, then there exists an associated unique global-in-time solution $\xi$ of \eqref{e:xis} such that $(\xi, \xi_s)\big|_{s=0} = (\xi_0, \xi_1)$, with $s\mapsto \cE_{\tot}^{\le J}(s)$ continuous and satisfying the decay estimate 
\begin{equation}
\label{e:maindecay}
\cE_{\tot}^{\le J}(s) \le C \cE^{\le J}_0 e^{-\kappa s}, \qquad \qquad s\ge 0.
\end{equation}
\end{THEOREM}

We prove this theorem in Section \ref{s:stability}, after we establish some preliminaries in Section \ref{s:Preliminaries} and prove the equivalence of energies $\cE_{\tot}^{\le J} \sim \cE_{\temp}^{\le J}$ mentioned above, in Section \ref{s:EquivEnergies}. We briefly discuss some related literature and our methodology before turning to these items.

Static self-gravitating elastic bodies have been studied under a range of constitutive hypotheses, for example \cite{BeigSchmidt2003Static,CalogeroLeonori2015GroundStates,AC2020,ACL2022}. Jia et al. investigated the multiplicity of spherical self-gravitating elastic equilibria and their linearized spectral stability under radial perturbations \cite{JiaKodioChapmanGoriely2019}. For dynamical bodies with a free boundary, Andersson, Oliynyk, and Schmidt established local well-posedness for compatible data near a relaxed configuration under a small-gravity assumption \cite{Andersson_2011}  based on an adaptation of Koch's well-posedness framework for mixed problems for nonlinear hyperbolic equations \cite{Koch1993}. Related local theory for nonlinear elastodynamic boundary problems can be traced back to \cite{DafermosHrusa1985,SHIBATA1989154}.  Elastic balls have been studied in the relativistic setting in \cite{ANPR2024}.

Explicit global motions of finite elastic bodies were recently constructed by Sideris, who obtained expanding and collapsing spherical solutions without gravity \cite{Sideris2023}, and by Calogero, who proved the existence of homologously expanding self-gravitating polytropic elastic balls \cite{Calogero2022}. Calogero's collapsing examples are supported by numerical evidence and have a different vacuum boundary condition. The authors constructed expanding, stationary, and collapsing separable self-gravitating hyperelastic bodies with a traction-free boundary \cite{JangLeslie2025}. In that setting, gravity can produce collapse even in the presence of positive residual pressure.

An important analog comes from the mass-critical Euler--Poisson system for gaseous stars. Homologous motions were investigated by Goldreich and Weber, Makino, and Fu and Lin \cite{GoWe1980,Makino1992,FuLin1998}; Hadžić and Jang later proved nonlinear radial stability of expanding stars \cite{HadzicJang2018CPAM}, with non-radial stability shown by Hadžić, Jang and Lam  \cite{HadzicJangLam2022Nonradial}.  For the compressible Euler equations without gravity, Sideris constructed expanding affine fluid motions \cite{Sideris2017Affine}, whose nonlinear stability was studied by Hadžić and Jang and by Shkoller and Sideris \cite{HadzicJang2018Inventiones,ShSi2019}. The gaseous vacuum boundary in these fluid problems differs from the traction-free solid boundary, where density remains positive. We also mention global results for elastic models on the whole space \cite{Sideris2000PrestressedWaves}. 

Our Theorem \ref{t:main} establishes global-in-time existence and stability for zero-energy perturbations of certain self-similarly expanding separable solutions of \eqref{e:Eulerian}, including a subclass of those separable solutions proved to exist by the authors in \cite{JangLeslie2025}. 
Following the framework developed in~\cite{HadzicJang2018CPAM}, we use coercivity of $\cL_\xi$ on mean-zero functions to control the mean-zero parts of time derivatives, the zero-energy constraint \mbox{\eqref{e:Energys}--\eqref{e:Js}} to control their means (cf. Lemma \ref{l:coercivity}), and elliptic estimates to recover spatial regularity (cf. Lemma \ref{l:EnEquiv}). The main structural differences come from the background profile and the nonlinear operator. While the enthalpy profile in~\cite{HadzicJang2018CPAM} satisfies a generalized Lane--Emden equation $w'' + \frac2R w' + \pi G w^3 = -\frac34 \mu, R\in [0,1]$ with $w'(0) = 0, \;w(1) = 0$, our profile equation~\eqref{e:separated} is strongly nonlinear, involves the ratio $y_* = R\varphi_R/\varphi$ (keeping in mind, however, that $w$ is a different kind of object than $\varphi$ is), and carries the nonlinear boundary condition $g'(y_*(1))=0$. Furthermore, our $\cL_\xi$ defined in~\eqref{e:Lmuxi} depends nonlinearly on $\xi$ and $\xi_R$ through the difference quotient $\frac{g'(y_*+z)-g'(y_*)}{z}$ inside the derivative. Its divergence structure nonetheless provides the crucial cancellations and coercivity needed to close the nonlinear energy estimates without loss of derivatives (cf. Remark \ref{r:Lform}). Lastly, we remark that our perturbation boundary conditions are simpler. While the vacuum boundary in~\cite{HadzicJang2018CPAM} motivates weighted Sobolev spaces adapted to the background, the traction-free conditions for both $\phi_*$ and $\phi=\phi_*(1+\xi)$ reduce to homogeneous Neumann conditions for $\xi$; see Section \ref{ss:Equations}. These conditions are essential at several points, particularly in the final integration by parts that closes the top-order estimates.

This paper deals exclusively with radial perturbations of self-similarly expanding separable solutions.  It would also be interesting to investigate non-radial stability, and to consider analogous questions for the other separable motions constructed in \cite{JangLeslie2025}. We leave these questions for future study.

\section{Preliminaries}

\label{s:Preliminaries}

In this section, we collect a number of technical preliminaries.  This will allow us to streamline our main proof in Section \ref{s:stability}, as well as the elliptic estimates in Section \ref{s:EquivEnergies}.  First, to simplify the technical difficulties of repeated spatial differentiation (which would otherwise be considerable, in light of the coordinate singularity $\frac1R$ that appears throughout our setup), Section \ref{ss:Preliminaries} develops the calculus of a class of operators that includes the $\cD^k$s; we use the notation $\cP_k$ for this class of operators. Section \ref{ss:Preliminaries}  also develops and records a number of basic estimates that will be used in what follows.  The framework here is inspired by that of \cite{GuoHadzicJang2020}, where a similar calculus was developed.  In order to reach our main proofs more quickly, we relegate to the Appendix the justification of any statements that take more than a few lines to establish.  Next, Section \ref{ss:backgroundsoln} shows that the background solution is sufficiently regular that the derivatives of the background quantities arising in the proof of \eqref{e:equivinduction} remain bounded under our assumptions and thus do not interfere with our estimates on $\xi$.  Finally, we briefly discuss the local well-posedness of our solutions to \eqref{e:xis} in Section \ref{ss:lwp}.

\subsection{Calculus for the class $\cP_k$}
\label{ss:Preliminaries}

In this subsection, we introduce and develop a calculus for the class $\cP_k$ of operators that we briefly mentioned above.  The latter includes the  $\cD^k$s defined in the Introduction and will play a prominent role in our proof of \eqref{e:equivinduction}.

When applying $\cD^k$ to \eqref{e:EnEquivStep1}, we will encounter the difficulty that the operators $\cD^k$ do not admit a natural product rule when $k\ge 3$.  This is a non-issue when $k=1$, since $\cD^1 = \partial_R$; furthermore, since $D_R(fg) = D_R f \cdot g + f \cdot \p_R g$, we have that 
\[
\cD^2(fg) = D_R \p_R f \cdot g + 2\p_R f \cdot \p_R g + f \cdot D_R \p_R g.
\]
However, the expansion of $\cD^3(fg)$ already cannot be written in terms of $\cD^j$, $j=0,1,2,3$ applied to $f$ and $g$ individually.  To correct this issue, we expand the class of operators under consideration and introduce the class $\cP_k$, defined as follows.  The first few are:
\[
\cP_0 = \{\id\}, 
\qquad 
\cP_1 = \{\p_R\},
\qquad 
\cP_2 = \{D_R \p_R, \;\tfrac1R \p_R\}.
\]
After $\cP_0$, they are defined recursively via
\[
\cP_{k+1} = 
\begin{cases}
	\{\p_R A : A\in \cP_{k}\}, & \text{ if } k \text{ is even}, \\
	\{D_R A : A\in \cP_{k}\} \cup \{\tfrac1R A : A\in \cP_{k}\}, & \text{ if } k \text{ is odd}. \\    
\end{cases}
\]

This expanded class of operators \textit{does} admit a natural product rule.  We record it as follows and provide a proof in the Appendix.  We do not track the values of the coefficients in this statement, as they are inconsequential for our purposes.

\begin{PROP}[Product Rule for operators in $\cP_k$] 
	\label{p:productrule}
	For $A\in \cP_k$, there exist constants $c^{ABC}$ such that 
	\begin{equation}
		\label{e:productrule}
		A(fg) = \sum_{\substack{\ell\in \{0,1,\ldots, k\}, \\ B\in \cP_\ell, \; C\in \cP_{k - \ell}}} c^{ABC} Bf \cdot Cg.
	\end{equation}   
\end{PROP}

Next, we record the commutators of elements of $\cP_k$ with $R\partial_R$, in the following Lemma.

\begin{LEMMA} The operators in $\cP_k$ admit commutators with $R \partial_R$ according to the following rules.
	\label{l:Dkcommutator}
    \label{l:Pkcommutator}
	\begin{enumerate}[label = (\alph*)]
	    \item We have 
	\begin{equation}
    \label{e:Dkcommutator}
		\cD^k \big(R \p_R \cdot\big)
		= (R\p_R + k) \cD^k = R \cD^{k+1} + \alpha_k \cD^k,
	\end{equation}
	where $\alpha_k = k$ if $k$ is even and $\alpha_k = k-4$ if $k$ is odd.
    \item If $A\in \cP_k\backslash \{\cD^k\}$, then 
	\[
	A(R\p_R \cdot) = \alpha A + B,
	\]
	where $\alpha\in \R$ and $B\in \cP_k$.  In particular, $A(R\p_R \cdot) \in \Span(\cP_k)$.
	\end{enumerate}
\end{LEMMA}

\begin{proof}
Statement (a) follows from the identities $\p_R\big( R \p_R \cdot \big) = (R\p_R + 1) \p_R$ and $D_R\big( R \p_R \cdot \big) = (R\p_R + 1) D_R$, together with an induction argument of just a few lines.  As for (b), when $A \ne \cD^k$, we can write $A = A'\frac1R \cD^j$ for some (odd) $j\in \{1, \ldots, k-1\}$ and some $A'\in \cP_{k-(j+1)}$.  Using \eqref{e:Dkcommutator}, 
	\[
	A(R\p_R \cdot) = A' \frac1R \cD^j (R\p_R \cdot) = A' \frac1R (R\cD^{j+1} + (j-4) \cD^j) = A' \cD^{j+1} + (j-4) A' \frac1R \cD^j = \alpha A + B,
	\]
	with $\alpha = j-4$ and $B = A' \cD^{j+1}$.
\end{proof}
The fact that the operators in $\cP_k$ admit a product rule is already, of course, quite useful.  The following Proposition tells us that, furthermore, working with $\cP_k$ rather than the $\cD^k$s alone costs almost nothing in both $L^2(R^4\dR)$ and $L^\infty$.  We give the proof in the Appendix; it relies on Lemma~\ref{l:Dkcommutator}.

\begin{PROP}
	\label{p:Pkcontrol} Assume $X\in C^k((0,1])$ and $A\in \cP_k$.
    \begin{enumerate}[label = (\alph*)]
    \item If $\cD^i X\in L^2(R^4\dR)$ for each $i\in \{0, \ldots, k\}$, then 
    \[
	\int_0^1 |AX|^2 R^4 \dR \lesssim \int_0^1 |\cD^k X|^2 R^4 \dR.
	\]
    \item If $\cD^i X\in L^\infty$ for each $i\in \{0, \ldots, k\}$, then 
    \[
\|AX\|_{L^\infty} \lesssim \|\cD^k X\|_{L^\infty}.
    \]
    \end{enumerate}
The implicit constants in both estimates depend only on $k$.
\end{PROP}

The following Proposition, which we prove in the Appendix, allows us to control the quantity $z = \frac{R\xi_R}{1+\xi}$ and its derivatives, which we will need to do repeatedly in Sections \ref{s:EquivEnergies} and \ref{s:stability}.  The upshot is that $L^2(R^4 \dR)$ estimates of $Az$ `cost' one more derivative than those of $A \xi$ (as might be expected, since $z$ already contains a factor of $\xi_R$), but $Az$ and $A \xi$ are on the same `perturbation energy level' from the perspective of $L^\infty$ estimates, due to the factor of $R$ in the numerator of $z = \frac{R\xi_R}{1+\xi}$ (which saves an application of~\eqref{e:Hardy} that would otherwise be needed to achieve a final weight of $R^4$).

\begin{PROP}
	\label{p:zderivs}
	Assume $A\in \cP_k$, with $k\ge 0$, and let $p$ be a nonnegative integer.  Assume $\xi = \xi(s,R)$ is regular enough so that all quantities below are finite, and also that  $1+\xi\ge \frac12$ on $[0,1]$.  Then for $z = \frac{R\xi_R}{1+\xi}$, and $\cE_{\tot}^{\le j}$ defined in terms of $\xi$ via \eqref{e:Etotdef}, we have 
		\begin{equation} 
			\label{e:Azest1}
			\int_0^1 |\partial_s^p A z|^2 R^4\dR \lesssim \cE_{\tot}^{\le p+k+1} \cdot \Big( 1 + \big( \cE_{\tot}^{\le p+k+1} \big)^{p+k} \,\Big).
		\end{equation} 
		Furthermore,
		\begin{equation} 
			\label{e:Azest2}
			\|\partial_s^p A z\|_{L^\infty}^2 \lesssim \cE_{\tot}^{\le p+k+3} \cdot \Big( 1 + \big(\cE_{\tot}^{\le p+k+2}\big)^{p+k} \,\Big).
		\end{equation} 
The implicit constants depend on $p$ and $k$ only.
\end{PROP}

The proof of Proposition \ref{p:zderivs} leverages the following Sobolev and Hardy inequalities, which are valid for sufficiently regular functions $X:(0,1]\to \mathbb{R}$. We record them here because they will also be used in what follows without further comment.
\begin{equation}
	\label{e:Sob}
	\|X\|_{L^\infty}^2 \lesssim \int_0^1 \big(|X|^2 + |\p_R X|^2\big) \dR,
\end{equation}
\begin{equation}
	\label{e:Hardy}
	\int_0^1 R^{2k} |X|^2 \dR \lesssim \int_0^1 R^{2(k+1)} (|X|^2 + |\p_R X|^2) \dR, \qquad k\ge 0.
\end{equation}

\begin{REMARK}
\label{r:R4}
The weight $R^4$ appears in most of the integrals throughout our analysis.  It is sometimes useful to view these weighted integrals instead as \textit{unweighted} integrals of radial functions on the unit ball in $\R^5$, which we denote by $B^5$.  More specifically, let $\xi$ belong to the space $\cH^k(R^4\dR)$, which we define via  
\begin{equation}
\label{e:Hkdef}
\|\xi\|_{\cH^k(R^4\dR)}^2:=\sum_{j=0}^k \|\cD^j \xi\|_{L^2(R^4\dR)}^2<+\infty.
\end{equation} 
Define $\Xi:B^5\to \R$ via $\Xi(x) = \xi(|x|)$.  Then $\|\Xi\|_{H^k(B^5)}$ (the Sobolev norm of $\Xi$ on $B^5$, defined in the usual way) is bounded above and below by (finite, positive) constant multiples of $\|\xi\|_{\cH^k(R^4\dR)}$. We make explicit use of this fact for $k=0,1,2$ in our estimates; it is obvious for $k=0,1$ and can be justified easily when $k=2$ using Proposition \ref{p:Pkcontrol}(a).  For larger $k$, it follows from elliptic regularity theory of the Laplace operator (note that $\cD^{2m}\xi(|x|) = \Delta^m \Xi(x)$).  In particular, this gives us access to the Sobolev embeddings $\cH^1(R^4\dR)\hookrightarrow L^{\frac{10}{3}}(R^4\dR)$ and $\cH^2(R^4\dR)\hookrightarrow L^{10}(R^4\dR)$, which we will use in Section~\ref{s:stability}.  One more fact that we will also use is that $\|\xi_R\|_{\cH^2(R^4\dR)}\lesssim \sum_{m=0}^3 \|\cD^m \xi\|_{L^2(R^4\dR)}$; this follows from the considerations discussed just above, together with the Hardy inequality \eqref{e:Hardy}.

On a related note, our nonlinear operator $\cL_\xi$ is, to leading order (and in the parameter regime of interest), essentially the five-dimensional radial Laplacian operator, and $\cE_{\tot}^{\le J}$ can be viewed as a sum of $H^j$ norms of radial functions on $B^5$.  This observation is what allows us to obtain local well-posedness of \eqref{e:xis} from existing literature rather than proving it directly.  See Section \ref{ss:lwp} for further discussion.
\end{REMARK}

\subsection{The background solution}
\label{ss:backgroundsoln}

As mentioned above, we will need to establish estimates on our background solution $\varphi$ in order to establish the necessary control over $\xi$ to prove \eqref{e:equivinduction}.  Theorem \ref{t:background} below provides the needed statement.  Before recording it, we pause to rewrite our equation for $\varphi$ in a form that will be more convenient below, and is valid whenever we can guarantee $g''>0$.  Following the derivation in \cite{JangLeslie2025}, equation \eqref{e:separated}$_1$ becomes 
\begin{equation}\label{e:background}
\frac{1}{R}\left(\p_R^2 \f + \frac{2}{R} \big(\p_R\f -\frac{\f}{R}\big) \right) = -\frac{1}{R^2}  \big(\p_R\f -\frac{\f}{R}\big) U (y_*) + V(\l, y_*),
\end{equation}
where 
\begin{equation} 
\label{e:Udef}
U(y_*) = 
\begin{cases} 
-(y_*+2)+\dfrac{g(y_*) + (4-y_*)g'(y_*)}{(y_*-1)g''(y_*)}, & y_*\ne 1, \\
0, & y_*=1,
\end{cases}
\end{equation}
\begin{equation} 
\label{e:Vdef}
V(\lambda, y_*) = \frac{\l}{\brho^{\frac13} g''(y_*)} \bigg( \frac{4\pi}{3} G \brho + \mu \lambda^3\bigg).
\end{equation} 
The functions $U(y_*)$ and $V(\lambda, y_*)$ are both as smooth in $y_*$ as $g$ allows (i.e., if $g\in C^{J+2}$, then $U$ and $V$ are $J$ times continuously differentiable in $y_*$ on the set where $g''\ne 0$), and $V(\lambda, y_*)$ is smooth in its first argument.  Equation \eqref{e:background} can be written in the more compact form
\begin{equation}\label{e:background1}
D_R\p_R \l  = -\frac{\p_R \l}{R}  U (y_*) + V(\l, y_*).  
\end{equation}

\begin{THEOREM}[Regularity of the Background Solution]
\label{t:background} 
Assume $g\in C^{J+2}((0,\infty))$, $J\in \N$, and that $g(1) = 1$, $g'(1) = -\frac13$, and $g''(1)>0$.  Let $\varphi\in C^2([0,1])$ be a solution of \eqref{e:separated} that satisfies 
	\begin{equation}
		\label{e:phiRRsmall}
		|\varphi''(R)| \le R \e,    \qquad R\in (0,1),
	\end{equation}
for some $\e>0$. 

\begin{enumerate}[label = (\alph*)]
    \item If $\e\in (0,1)$, then 
    	\begin{equation}
		\label{e:y*small}
		|\lambda-1| \le \frac13 \e, \qquad |y_* - 1|\le \frac12 \e,
        \qquad |y_*'| \le 3\e. 
	\end{equation}
    \item There exists $\e_1\in (0,1)$, depending only on $g$, such that if $\e\in (0,\e_1]$, then 
    \begin{equation}\label{e:regularity}
		\sum_{\substack{A\in \cP_{\ell},\\ \ell\in \{0, \ldots, J+2\}}} \|A \l \|_{L^\infty} + \sum_{\substack{B\in \cP_{\ell},\\ \ell\in \{0, \ldots, J+1\}}} \|B y_*\|_{L^\infty}<\infty.  
	\end{equation}
\end{enumerate}
\end{THEOREM}

\begin{REMARK}
\label{r:g''requirement}
The requirement \eqref{e:phiRRsmall} for a solution $\varphi$ of \eqref{e:separated} is known to be justified under more stringent assumptions on $g$ (and $\mu$).	In particular, the authors have shown in \cite{JangLeslie2025} that, under the assumption 
\begin{equation}
\label{e:sepexisthyp}
g''(1) \ge \max\bigg\{ 50\Big( 50 + \sup_{|y-1|\le \frac12} |g'''(y)|\Big), \, \frac{10}{\e}\bigg\},		
\end{equation}
there exists $\mu_0 = \mu_0(\e)>0$ such that whenever $|\mu|<\mu_0$, there exists $\brho>0$ and a corresponding classical solution $\varphi$ of \eqref{e:separated} (equipped with the parameters $\mu$ and $\brho$) that satisfies \eqref{e:phiRRsmall}.  However, the assumption \eqref{e:sepexisthyp} is not optimal (and probably far from it), which is why we do not build it directly into the statement of Theorem \ref{t:background}.  
\end{REMARK}

\begin{proof}[Proof of Theorem \ref{t:background}] We begin with the proof of statement (a).  Assume $\e\in (0,1)$.  
Since $\varphi(0) = \varphi''(0) = 0$ and $\varphi(1) = 1$, we may write 
\begin{equation}
\varphi(R) = R \varphi'(0) + \int_0^R (R-\tau)\varphi''(\tau) \dtau,
\qquad \varphi'(0) = 1 - \int_0^1 (1 - \tau) \varphi''(\tau) \dtau.
\end{equation}
Some elementary manipulations and an application of \eqref{e:phiRRsmall} then yield 
\begin{equation}
\label{e:ystarminus1est}
|y_*(R)-1|\le \frac12 R^2\e,
\end{equation}
\begin{equation}
\label{e:lambdaprimeest}
|\lambda(R)-1|\le \frac13 \e,
\qquad 
|\lambda'(R)|\le \frac13 R\e. 
\end{equation}
Combining \eqref{e:ystarminus1est} and \eqref{e:lambdaprimeest} allows us to conclude that  $|y_*'(R)| = \big| \frac{\varphi''(R) - y_*(R) \lambda'(R)}{\lambda(R)}\big| \le  3 R\e$, so that in particular, \eqref{e:y*small} holds.  Statement (a) is proved.

As for statement (b), let us fix $\e_1\in (0,1)$ small enough that $g''(y) \ge \frac12 g''(1)$ for $y\in (1-\e_1, 1+\e_1)$.  We prove \eqref{e:regularity} by induction.  More specifically, we show that for any $m\in \{0, \ldots, J\}$, we have $A\lambda\in L^\infty$ whenever $A\in \cP_{\ell}$, $\ell\le m+2$.  Note that, since $y_* - 1 = (R\partial_R \lambda) \lambda^{-1}$, the inductive hypothesis on $\lambda$ (together with Lemma \ref{l:Dkcommutator}) implies that $By_*\in L^\infty$ whenever $B\in \cP_\ell$, $\ell\le m+1$; therefore, we need only run our induction on the regularity of $\lambda$.  Moreover, equation \eqref{e:background1}, together with the regularity of $U$ and $V$, implies that $\lambda \in C^{J+2}((0,1])$, so that we only need to worry about the boundedness of the $\cD^j \lambda$'s near $R=0$.  Since we have already established above that $|\lambda'(R)|\lesssim R \e$, equation \eqref{e:background1} also allows us to conclude that $D_R \partial_R \lambda$ is bounded, so that our inductive hypothesis holds for $m = 0$.

Next, we assume the desired conclusion holds for $m=K-1$ (with $K\in \{1, \ldots, J\}$) and apply $\cD^K$ to \eqref{e:background1}.  The result can be written 
\begin{equation}
\cD^{K+2} \lambda = -\frac{U(y_*)}{R} \bigg[ R\cD^K \bigg( \frac1R \lambda_R \bigg)\bigg] + \cQ_K, 
\end{equation} 
where 
\begin{equation}
\cQ_K = -[\cD^K, U(y_*)] \bigg( \frac1R \lambda_R \bigg) + \cD^K \big( V(\lambda, y_*) \big),
\end{equation}
and $[\cD^K, U(y_*)]$ denotes the commutator $\cD^K \big( U(y_*)\cdot\big) - U(y_*) \cD^K$.  The quantity $\cQ_K$ is bounded, by our inductive hypothesis.  The quantity $\frac{U(y_*(R))}{R}$ is also bounded. (In fact, it is even $O(R)$, since $U$ is $C^1$ near $1$, with $U(1) = 0$, so that $|U(y_*(R))|\lesssim R^2 \e$ by \eqref{e:ystarminus1est}.)  Furthermore, we have
\[
R\cD^K \bigg( \frac1R \lambda_R \bigg)
= \cD^{K+1}\lambda -(5+\alpha_{K-1}) \cD^{K-1} \bigg( \frac1R \partial_R \lambda \bigg),
\]
where $\alpha_{K-1}$ is either $K-1$ or $K-5$, as in Lemma \ref{l:Dkcommutator}(a).  Since this quantity is also bounded, we conclude that $\cD^{K+2}\lambda$ is bounded.  We invoke Proposition \ref{p:Pkcontrol}(b) to complete the proof.
\end{proof}

\subsection{Local-in-time existence of solutions}
\label{ss:lwp}
We now briefly discuss the local-in-time well-posedness of the system \eqref{e:xis}.  The key thing to note here is that, in the parameter regime of interest, $\cL_\xi$ is essentially the five-dimensional radial Laplacian at leading order, and we can therefore view it as a nonlinear wave equation.  Theorem \ref{t:lwp} below is then a consequence of the work \cite{SHIBATA1989154} of  Shibata and Kikuchi (more specifically, Theorem 8.1 in that paper), which treats the local well-posedness of nonlinear scalar wave equations with Neumann-type boundary conditions. 

\begin{THEOREM}
\label{t:lwp}
Assume $g\in C^{J+2}((0,\infty))$, $J\ge 5$ and that $g(1) = 1$, $g'(1) = -\frac13$, and $g''(1)>0$.  There exists $\e_L>0$, depending only on $g$ and $J$, such that the following statement holds whenever $\varphi$ is a classical solution of \eqref{e:separated} (with some associated parameters $\mu<0$ and $\brho>0$) satisfying  $|\varphi''(R)|\le R\e_L$,  $R\in (0,1)$.

There exist $\th_L,T_L,C_L>0$ (depending on $g$, $J$, $\varphi$, $\brho$, and $\mu$) with the following property. Given $(\xi_0, \xi_1)$ which satisfies $\cE_0^{\le J}\le \theta_L$ and is compatible with the boundary conditions of the system in the sense described in Section \ref{sss:IC}, there exists a unique solution $\xi$ solving \eqref{e:xis} on the time interval $[0,T_L]$ and satisfying the initial conditions $(\xi, \xi_s)\big|_{s=0} = (\xi_0, \xi_1)$.  Furthermore, the function $s\mapsto \cE_{\tot}^{\le J}(s)$ is continuous on $[0,T_L]$ for this solution, and it satisfies the energy bound
\begin{equation} 
\label{e:lwpenergybd}
\sup_{\tau\in[0,T_L]}\cE_{\tot}^{\le J}(\tau)\le C_L \cE_{0}^{\le J}.
\end{equation} 
\end{THEOREM}

\section{Equivalence of energies}

\label{s:EquivEnergies}

In this section, we prove the following lemma, which states (for $J\ge 5$) that $\cE^{\le J}_{\temp}$ and $\cE^{\le J}_{\tot}$ are comparable when the latter is small enough and when $\xi$ solves \eqref{e:xis}.  As mentioned already in the introduction, this will allow us to replace $\cE_{\temp}^{\le J}$ with $\cE_{\tot}^{\le J}$ in a crucial place in the proof of Theorem \ref{t:Etotest}.

\begin{LEMMA}
\label{l:EnEquiv}
Assume $g\in C^{J+2}((0,\infty))$, $J\ge 5$ and that $g(1) = 1$, $g'(1) = -\frac13$, and $g''(1)>0$.  There exists $\e_2>0$, depending only on $g$ and $J$, such that the following statement holds whenever $\varphi$ is a classical solution of \eqref{e:separated} (with some associated parameters $\mu<0$ and $\brho>0$) satisfying  $|\varphi''(R)|\le R\e_2$,  $R\in (0,1)$.

There exist positive constants $\th_2$ and $C_2$ (depending on $\e_2$, $g$, $J$, $\varphi$, $\brho$, and $\mu$), such that if $\xi$ is a solution of \eqref{e:xis}, then at any time $s$ for which $\cE^{\le J}_{\tot}(s) \le \th_2$, it follows that $\cE^{\le J}_{\tot}(s) \le C_2 \,\cE^{\le J}_{\temp}(s)$. 
\end{LEMMA}

We briefly comment on our strategy before diving into the proof.  We will show something slightly stronger than what is claimed, namely that, whenever $g$ is $J+2$ times continuously differentiable, with $J\ge \max\{ j+k, \lfloor \frac{j+k}{2} \rfloor + 3\}$, we have 
\begin{equation}
\label{e:equivinduction}
\|\partial_s^j \cD^k \xi\|_{L^2(R^4 \dR)}^2 \lesssim \cE_{\temp}^{\le j+k} + \cE_{\tot}^{\le j+k-1} \cE_{\tot}^{\le \lfloor \frac{j+k}{2} \rfloor + 3}.
\end{equation}
Whenever $\lfloor \frac{J}{2} \rfloor + 3 \le J$, i.e., whenever $J \ge 5$, the inequality \eqref{e:equivinduction} implies that 
\[
\cE_{\tot}^{\le J} \lesssim \cE_{\temp}^{\le J} + \cE_{\tot}^{\le J-1} \cE_{\tot}^{\le J},
\]
from which the conclusion of Lemma \ref{l:EnEquiv} follows, upon taking $\th_2$ small enough to absorb the quadratic term into the left side.  

The cases $k = 0,1$ of \eqref{e:equivinduction} are trivial to verify, since then $\|\partial_s^j \cD^k \xi\|_{L^2(R^4 \dR)}^2$ is one of the terms in the definition of $\cE_{\temp}^{\le j+k}$.  (We interpret $\cE_{\tot}^{\le -1}$ as zero when it arises.)
To prove \eqref{e:equivinduction} when $k\ge 2$, we differentiate \eqref{e:xis}$_1$ repeatedly and induct on the number of spatial derivatives.  The version of \eqref{e:xis}$_1$ that we use for this purpose is 
\begin{equation} 
	\label{e:EnEquivStep1}
	F_1 D_R \xi_R = \xi_{ss} - \frac12 b\xi_s + 3 \mu \xi - F_2 \cdot \frac{\xi_R}{R} - F_3 \xi^2 - F_4 \xi_R^2,
\end{equation}
where (writing $\d_z h(y_*) = \frac{h(y_* + z)- h(y_*)}{z}$) the functions $F_j$, $j=1,2,3,4$, are defined as follows.   
\begin{equation} 
\label{e:F1F2def}
F_1 = \frac{\brho^{\frac13}}{\l^3(1+\xi)^3} g''(y_* + z), 
\qquad 
F_2 = \frac{\brho^{\frac13}}{\l^3(1+\xi)^3} \big( \d_z g''(y_*)(R y_*') - \d_z g'(y_*) (y_*-1) \big), 
\end{equation} 
\begin{equation}
\label{e:F3F4def}
F_3 = \cN_1(\xi),
\qquad 
F_4 = -\frac{\brho^{\frac13}}{\l^3(1+\xi)^4}\Big[ g''(y_*+z) + \p_z\d_z g(y_*) + 3\p_z\d_z g'(y_*)\Big]
\end{equation}

Equation \eqref{e:EnEquivStep1} is obtained from expanding the definition \eqref{e:Lmuxi} of $\cL_{\xi}$ in \eqref{e:xis}$_1$.  The carefully chosen form \eqref{e:Lmuxi} of $\cL_{\xi}$ plays a key role in the proof of Theorem \ref{t:Etotest}, but its expanded form is more useful for the proof of Lemma \ref{l:EnEquiv}. 

The main strategy of our estimates leading to \eqref{e:equivinduction} can already be understood from the form of \eqref{e:EnEquivStep1}, with a bit of additional explanation about the roles of each of the functions $F_j$.  At the $k$th step of the induction, we apply $\partial_s^j \cD^{k-2}$ to \eqref{e:EnEquivStep1}, $j=0,1,\ldots, J-k$.  When $y_*-1$, $Ry_*'$, $\xi$, and $z$ are small enough, we have that $\underline{F_1}:= \inf F_1 > 0$ and that $|F_2|\ll \underline{F_1}$.  Because of the coordinate singularity, $\frac{\xi_R}{R}$ and $D_R \xi_R$ have effectively the same derivative count, and therefore, the two contributors to the left side of~\eqref{e:equivinduction} are $F_1\partial_s^j \cD^{k}\xi$ and $F_2 \partial_s^j \cD^{k-2}\big( \frac{\xi_R}{R} \big)$.  The remaining terms in the expansions of $\partial_s^j \cD^{k-2}\big( F_1 D_R \xi_R\big)$ and \mbox{$\partial_s^j \cD^{k-2}\big( F_2 \cdot\frac{\xi_R}{R}\big)$} are at least quadratic in $\xi$ and its (sufficiently low-order) derivatives, and they contribute, along with $\partial_s^j \cD^{k-2}(F_3 \xi^2)$ and $\partial_s^j \cD^{k-2}(F_4 \xi_R^2)$, toward the $\cE_{\tot}^{\le j+k-1} \cE_{\tot}^{\le \lfloor \frac{j+k}{2} \rfloor + 3}$ in \eqref{e:equivinduction}.  Of course, the appearance of $\cE^{\le j+k}_{\temp}$ in \eqref{e:equivinduction} is due to the remaining linear terms $\partial_s^j \cD^{k-2}(\xi_{ss} - \frac12 b \xi_s + 3\mu \xi)$.

\begin{proof}[Proof of Lemma \ref{l:EnEquiv}]

We are ready to proceed with our estimates on $\partial_s^j \cD^k \xi$ and prove the estimate \eqref{e:equivinduction}.   The case $k=2$ is slightly different from $k\ge 3$, so we treat the former individually before moving on.  

\underline{\textbf{Estimates for $k=2$}}.
We first estimate \eqref{e:EnEquivStep1} directly:
\begin{align*}
	\int_0^1 |D_R \xi_R|^2 R^4 \dR 
	& \lesssim \int_0^1 \big( |\xi_{ss}|^2 + |\xi_s|^2 + |\xi|^2 \big) R^4 \dR + \|F_2\|_{L^\infty}^2 \int_0^1 |\cD^2 \xi|^2 R^4 \dR \\
	& \qquad \quad + \|\xi\|_{L^\infty}^2 \int_0^1 |\xi|^2 R^4 \dR + \|\xi_R\|_{L^\infty}^2 \int_0^1 |\xi_R|^2 R^4 \dR 
\end{align*}
Theorem \ref{t:background} (together with the definition of $F_2$ in \eqref{e:F1F2def}) guarantees that by choosing $\e_2>0$ small enough, we can ensure that $\|F_2\|_{L^\infty}^2$ is very small compared to the implied constant in this inequality, so that we can absorb this term into the left side of the equation.  In some future steps, we will perform this manipulation (possibly after shrinking $\e_2$ again, a total of finitely many times) without further comment. 

The estimate becomes 
\[
\int_0^1 |\cD^2 \xi|^2 R^4 \dR \lesssim \cE_{\temp}^{\le 2} + \cE_{\tot}^{\le 1}\cE_{\tot}^{\le 4},
\]
which is exactly what is required for \eqref{e:equivinduction} when $j=0$ and $k=2$.  

Now, applying time derivatives, we may write the following (for some constants $a_0, \ldots, a_j$):
\[
\begin{split}
	F_1 \p_s^j \cD^2 \xi 
	& = \p_s^j \bigg(\xi_{ss} - \frac12 b\xi_s + 3\mu \xi \bigg) - F_2 \p_s^j \bigg( \frac{\xi_R}{R} \bigg) \\
	& \qquad + \sum_{\ell=1}^j a_\ell \bigg[ (\p_s^\ell F_1)(\p_s^{j-\ell} \cD^2 \xi) + (\p_s^\ell F_2)\bigg( \frac{\p_s^{j-\ell}\xi_R}{R} \bigg) \bigg] \\
	& \qquad + \sum_{\ell=0}^j a_\ell \big[ (\p_s^\ell F_3)(\p_s^{j-\ell} (\xi^2)) + (\p_s^\ell F_4)(\p_s^{j-\ell} (\xi_R^2)) \big] \\
\end{split}
\]
When we take $L^2(R^4\dR)$ norms, the linear $F_2$ term in the first line can be absorbed into the left side as before; the squares of the remaining terms on the right side of the equation are bounded above by $\cE^{\le j+2}_{\temp}$.  The terms in the second and third lines are quadratic.  We explicitly consider only the terms of the form $(\p_s^\ell F_1)(\p_s^{j-\ell} \cD^2 \xi)$, $1\le \ell \le j$.  (Once we treat these, the other term in the second line can be dealt with similarly, and the terms in the third line can also be treated in a similar fashion but more easily.)  

In the term $(\p_s^\ell F_1)(\p_s^{j-\ell} \cD^2 \xi)$, $1\le \ell \le j$, the total number of derivatives that can fall on a single $\xi$ is $j+1$, if we count the two spatial derivatives $\cD^2$. (This also counts $\p_s^j z$ as a term with a total of $j+1$ derivatives on $\xi$.) The smaller of $\ell$ and $j-\ell + 2$ cannot exceed $\lfloor \frac{j+2}{2} \rfloor$, so the `worst case' scenario for the second-highest total number of derivatives is the one where $\lfloor \frac{j+2}{2} \rfloor$ derivatives in $s$ fall on either $\xi$ or $z$.  Putting an $L^2(R^4\dR)$ norm on the highest derivative and an $L^\infty$ norm on the second-highest one, we get a total contribution of no worse than $\cE^{\le j+1}_{\tot} \cE^{\le \lfloor \frac{j+2}{2} \rfloor + 3}_{\tot}$ from these terms.  Note that our smallness assumption allows us to restrict attention to consideration of the two highest-order derivatives involved in any term which is quadratic or higher-order in powers of $\xi$, $z$, and their derivatives.  Therefore, for $0\le j \le J-2$, we have
\begin{equation}
	\label{e:sderivest1}
	\int_0^1 |\p_s^j \cD^2 \xi|^2 R^4 \dR \lesssim \cE_{\temp}^{\le j+2} + \cE^{\le j+1}_{\tot} \cE^{\le \lfloor \frac{j+2}{2} \rfloor + 3}_{\tot},
\end{equation}
as needed.  

\underline{\textbf{Estimates for $k\ge 3$}}.
Now, suppose the desired statement holds for $k=0,1,\ldots, K-1$; we show that it must be true when $k = K$, with $K\ge 3$.  To this end, we apply $\partial_s^j \cD^{K-2}$ to \eqref{e:EnEquivStep1} to obtain the following:
\begin{equation}
	\label{e:DK-2}
	\begin{split}
		F_1 \partial_s^j \cD^{K} \xi 
		& = \partial_s^j \cD^{K-2} \big( \xi_{ss} - \frac12 b \xi_s + 3 \mu \xi \big) - F_2 \partial_s^j \cD^{K-2} \bigg( \frac{\xi_R}{R} \bigg) - \partial_s^j \cD^{K-2}\big( F_3 \xi^2 + F_4 \xi_R^2 \big) \\
        & - [\partial_s^j \cD^{K-2}, F_1](D_R \xi_R) - [\partial_s^j \cD^{K-2}, F_2]\bigg( \frac{\xi_R}{R} \bigg),
	\end{split}
\end{equation}
where $[\partial_s^j \cD^{K-2}, F_1]$ is the commutator of $\partial_s^j \cD^{K-2}$ with $F_1$, and similarly for $[\partial_s^j \cD^{K-2}, F_2]$.  The terms on the right side of the first line of \eqref{e:DK-2} are estimated via Proposition \ref{p:Pkcontrol} and our inductive hypothesis, as well as \eqref{e:y*small} and the resulting bound $\|F_2\|_{L^\infty}\lesssim \e_2$.  More specifically, 
\begin{align*}
	& \int_0^1 \big( |\partial_s^{j+2} \cD^{K-2} \xi|^2 +  |\partial_s^{j+1} \cD^{K-2} \xi|^2 + |\partial_s^{j} \cD^{K-2} \xi|^2\big) \,R^4\dR +  \bigg\| F_2 \partial_s^j \cD^{K-2} \bigg( \frac{\xi_R}{R} \bigg) \bigg\|_{L^2(R^4 \dR)}^2 \\
	& \lesssim \cE_{\temp}^{\le j+K} + \cE_{\tot}^{\le j+K-1}\cE_{\tot}^{\le \lfloor \frac{j+K}{2} \rfloor + 3} + \e_2^2 \|\partial_s^j \cD^K \xi\|_{L^2(R^4 \dR)}^2.
\end{align*}
Decreasing $\e_2$ if necessary, we may absorb $\e_2^2 \|\partial_s^j \cD^K \xi\|_{L^2(R^4 \dR)}^2$ into the left side of our overall estimate.  

We now turn to the commutator term $[\partial_s^j \cD^{K-2}, F_1](D_R \xi_R)$, the expansion of which yields contributions of the form $(\partial_s^m BF_1)(\partial_s^{j-m} CD_R \xi_R)$, where $B\in \cP_\ell$, $C\in \cP_{K-2-\ell}$, and $0\le m\le j$, $0\le \ell\le K-2$, and $m+\ell\ge 1$.  (The other commutator can be dealt with in an entirely similar way; we leave the details to the reader.) Expanding $\partial_s^m BF_1$ in the expression at hand gives rise to a continuous function of $\lambda$, $y_*$, $Ry_*'$, $\xi$, and $z$ (namely, some partial derivative of $F_1$ with respect to its arguments), multiplied by powers of quantities like $B_1 \lambda$, $B_2 y_*$, $B_3 (Ry_*')$, $\partial_s^{\ell_4} B_4 \xi$, and $\partial_s^{\ell_5} B_5 z$.  The only way for the corresponding terms (including their contributions from $CD_R \xi_R$) to be linear in $\xi$ and its derivatives is for $m$ to be equal to zero and for all of the (at most $\ell$) spatial derivatives from $B$ to ultimately fall on the background quantities $\lambda$, $y_*$, $Ry_*'$.  By Theorem \ref{t:background}, these linear terms have a contribution of $\lesssim \sum_{p=0}^{K-1} \|\partial_s^{j} \cD^{p} \xi\|_{L^2(R^4\dR)}^2$, and our induction hypothesis guarantees that this quantity satisfies the desired bound.  

As for the \textit{nonlinear} terms of the form  $(\partial_s^m BF_1)(\partial_s^{j-m} CD_R \xi_R)$, we note that the highest number of (time or spatial) derivatives that can land on $\xi$ is $j+K-1$; the second-largest number of derivatives cannot exceed $\lfloor \frac{j+K}{2} \rfloor$.  We put the $L^2(R^4 \dR)$ norm on the contribution with the highest number of derivatives and achieve a bound of the form $\cE_{\tot}^{\le j+K-1} \cE_{\tot}^{\le \lfloor \frac{j+K}{2} \rfloor + 3}$ on all such terms.  We conclude that the estimate \eqref{e:equivinduction} holds for $k=K$, which completes the proof of Lemma \ref{l:EnEquiv}.  
\end{proof}

\section{Stability}

\label{s:stability}

In this section, we prove our main existence and stability result, Theorem \ref{t:main}.  Most of the remaining work is in establishing the following higher-order energy estimate.
\begin{THEOREM}
\label{t:Etotest}
Assume $g\in C^{J+2}((0,\infty))$, $J\ge 5$, and that $g(1) = 1$, $g'(1) = -\frac13$, and $g''(1)>0$.  There exists $\e_0>0$ depending only on $g$ and $J$, and there exists an absolute constant $c_0>0$, such that the following holds.  

Given $\brho>0$ and $\mu\in (-c_0\brho^{\frac13}g''(1),0)$, let $\varphi$ be a classical solution of \eqref{e:separated} (equipped with the parameters $\brho$, $\mu$) that satisfies  $|\varphi''(R)|\le R\e_0$, $R\in (0,1)$.  Then there exist positive constants $\theta_0$, $\kappa_0$, and $\kappa_1$, depending on $g$, $J$, $\varphi$, $\mu$, and $\brho$, such that if $\xi$ is a solution of \eqref{e:xis} on some time interval $[0,T]$, satisfying the zero-energy condition \eqref{e:Energys}--\eqref{e:Js} and the a priori estimate  
\[
\sup_{\tau\in[0,T]}\cE_{\tot}^{\le J}(\tau)\le \theta_0,
\]
then for all $0\le \widetilde{s}<s\le T$, we have
\begin{equation} 
\label{e:Etotdecay1}
\cE_{\tot}^{\le J}(s)+\kappa_0\int_{\widetilde{s}}^s\cE_{\tot}^{\le J}(\tau)\,\dtau
\ \le\ \kappa_1\,\cE_{\tot}^{\le J}(\widetilde{s}).
\end{equation} 
\end{THEOREM}

Assume provisionally that Theorem \ref{t:Etotest} holds.  We prove our main theorem as an easy consequence of its conclusion.

\begin{proof}[Proof of Theorem \ref{t:main}]

Let $\e_0$, $c_0$, $\theta_0$, $\kappa_0$, and $\kappa_1$ be as in Theorem \ref{t:Etotest};  let $\e_L$, $\theta_L$, $T_L$, and $C_L$ be as in Theorem \ref{t:lwp}.  Fix $\e = \min\{\e_0, \e_L\}$, $c=c_0$, and $\theta = \min\{\frac12\theta_0, \theta_L\} / \big( \max\{C_L, 1\}\cdot \max\{\kappa_1, 1\}\big)$.  Let $(\xi_0, \xi_1)$ be initial data which is compatible with the boundary conditions of \eqref{e:xis} and satisfies the zero-energy condition \eqref{e:Energys}--\eqref{e:Js} and the smallness condition $\cE_0^{\le J} \le \theta$; let $\xi$ be the unique solution of \eqref{e:xis} generated from $(\xi_0, \xi_1)$, which is guaranteed to exist at least for a short time by Theorem \ref{t:lwp}. Let $T_{\max}$ be the supremum of all times $T$ such that $\xi$ exists as a solution of \eqref{e:xis} on $[0,T]$ and satisfies $\sup_{\tau\in [0,T]} \cE^{\le J}_{\tot}(\tau)\le \theta_0$.  By the energy bound \eqref{e:lwpenergybd} of Theorem \ref{t:lwp}, and the fact that $C_L\theta\le \theta_0$, we have $T_{\max}\ge T_L$.  In fact, by the same argument, applied to the configuration $(\xi(s), \xi_s(s))$ rather than $(\xi_0, \xi_1)$ at any time $s\in (0,T_{\max})$, we may conclude that $T_{\max} \ge s + T_L$.  (Indeed, $(\xi(s), \xi_s(s))$ still satisfies the compatibility conditions and the zero-energy requirement; furthermore, by our choice of $\theta$, we have $\cE_{\tot}^{\le J}(s) \le \theta_L$, so that the solution exists on $[s,s+T_L]$, and also $C_L \cE_{\tot}(s)\le \frac12\theta_0$, so that $\cE_{\tot}^{\le J} \le \theta_0$ on $[s, s+T_L]$.) We can therefore bootstrap indefinitely to conclude that $T_{\max} = +\infty$, which establishes global-in-time existence of $\xi$ and the bound $\cE_{\tot}^{\le J} \le \theta_0$.

As for the decay estimate \eqref{e:maindecay}, we note that \eqref{e:Etotdecay1} implies that $\cE_{\tot}^{\le J}(s)\le \kappa_1 \cE_{\tot}^{\le J}(\tau)$ for all $\tau\in [0,s]$.  Integrating then yields
\[
\cE^{\le J}_{\tot}(s) \le \frac{\kappa_1}{s-\widetilde{s}} \int_{\widetilde{s}}^s \cE_{\tot}^{\le J}(\tau) \dtau \le \frac{\kappa_1}{s-\widetilde{s}} \cdot \frac{\kappa_1}{\kappa_0} \cE_{\tot}^{\le J}(\widetilde{s}),
\qquad 0\le \widetilde{s}<s.
\]
Set $T_0 = \frac{2\kappa_1^2}{\kappa_0}$.  Then for $s\in ((\ell-1) T_0, \ell T_0]$, $\ell\in \N$, this last inequality gives us 
\begin{align*}
\cE_{\tot}^{\le J}(s) \le \kappa_1 \cE_{\tot}^{\le J}((\ell-1) T_0) \le \frac12 \kappa_1 \cE_{\tot}^{\le J}((\ell-2) T_0)\le \cdots \le 2^{-(\ell-1)} \kappa_1 \cE_0^{\le J},  
\end{align*}
that is, 
\[
\cE_{\tot}^{\le J}(s) \le 2\kappa_1 e^{-\frac{\log 2}{T_0}s} \cE_0^{\le J},
\qquad s\ge 0.
\]
This is exactly \eqref{e:maindecay}, with $C = 2\kappa_1$ and $\kappa = \frac{\log 2}{T_0} = \frac{\kappa_0}{2\kappa_1^2} \log 2$.
\end{proof} 

The rest of the section is dedicated to the proof of Theorem \ref{t:Etotest}.  

\subsection{Higher-order energy estimate}

\label{ss:higherorderenergyestimate}

Our starting point in the proof of Theorem \ref{t:Etotest} is the following time-differentiated version of \eqref{e:xis}, writing $\xi^{(j)} := \partial_s^j \xi$ and using $[\partial_s^j, \cL_{\xi}] = \partial_s^j \cL_{\xi} - \cL_{\xi} \partial_s^j$ to denote the commutator of $\partial_s^j$ and $\cL_{\xi}$:
\begin{equation}
\label{e:xistimederiv}
\xi^{(j)}_{ss} - \frac12 b \xi^{(j)}_s + (\cL_{\xi} + 3 \mu)\xi^{(j)} = \partial_s^j \cN(\xi) - [\partial_s^j, \cL_{\xi}](\xi).
\end{equation}
Let us leverage \eqref{e:xistimederiv} to obtain a higher-order energy identity involving weighted $L^2$  norms of the derivatives of $\xi$.  We use the weight $R^4 \lambda^2 (1+\xi)$, which is chosen so that $\cL_{\xi}$ is self-adjoint with respect to the associated $L^2$ inner product (provided that $\xi$ satisfies the Neumann boundary conditions \eqref{e:xis}$_2$).  We set notation accordingly:
\[
\langle u,v\rangle_\xi = \int_0^1 u(R) v(R) R^4 \lambda^2 (1+\xi)\dR,
\qquad \|u\|_\xi^2 = \langle u, u\rangle_\xi.
\]
In the regime of interest, we have $|\xi|\le \frac12$ and $|\lambda-1|\le \frac13$, so that the norm $\|\cdot\|_\xi$ is bounded above and below by absolute constant multiples of $\|\cdot\|_{L^2(R^4 \dR)}$.

We take the $\langle\cdot, \cdot \rangle_\xi$ inner product of \eqref{e:xistimederiv} with $2\xi^{(j+1)} - b k_0 \xi^{(j)}$, where $k_0$ is a small, positive constant, the precise value of which is determined below.  After an additional integration, the result can be written as 
\begin{equation} 
\label{e:xisjIP}
\cF_{j,1}(\tau)\bigg|^{\tau = s}_{\tau = \widetilde{s}} + |b| \int_{\widetilde{s}}^s \cF_{j,2}(\tau)\,\dtau = \int_{\widetilde{s}}^s \cR_j(\tau) \,\dtau,
\end{equation}
where 
\[
\cF_{j,1} = \|\xi^{(j+1)}\|_\xi^2 + \frac{b^2 k_0}{4} \|\xi^{(j)}\|_\xi^2 - b k_0 \langle \xi^{(j+1)}, \xi^{(j)}\rangle_\xi + \langle (\cL_{\xi} + 3\mu)(\xi^{(j)}), \xi^{(j)} \rangle_\xi,  
\]
\[
\cF_{j,2} = (1-k_0) \|\xi^{(j+1)}\|_\xi^2 + k_0 \langle (\cL_{\xi} + 3 \mu)(\xi^{(j)}), \xi^{(j)}\rangle_\xi,
\]
\begin{equation} 
\label{e:Rj}
\begin{split} 
\cR_j & = \big\langle \partial_s^j \cN(\xi) - [\partial_s^j, \cL_{\xi}](\xi), 2\xi^{(j+1)} - b k_0 \xi^{(j)} \big\rangle_\xi 
+ \brho^{\frac13} \int_0^1 \partial_s \bigg( \frac{\delta_z g'(y_*)}{\lambda(1+\xi)^2} \bigg) |\xi^{(j)}_R|^2 R^4 \dR \\
& \qquad +  \int_0^1 \bigg [ |\xi^{(j+1)}|^2 - \frac{b^2}{4} ( 6 - k_0 ) |\xi^{(j)}|^2 - b k_0 \xi^{(j+1)} \xi^{(j)} \bigg] \lambda^2 R^4 \xi_s \dR.
\end{split} 
\end{equation} 
Note that $\cF_{j,1}$ and $\cF_{j,2}$ satisfy the obvious bound  
\begin{equation}
\label{e:Fjellupper}
\cF_{j,\ell}\lesssim \cE^{\le j+1}_{\temp}, \qquad \ell = 1,2.
\end{equation} 

We are now in a position to describe the remaining elements of our proof.  With Lemma \ref{l:EnEquiv} and the identity \eqref{e:xistimederiv} in hand, there are essentially two remaining technical steps.  The first is to establish coercivity estimates on $\cF_{j,1}$ and $\cF_{j,2}$, of the form
\begin{equation}
\label{e:Fjellcoercive}
\cF_{j,\ell} \ge c\big(\|\partial_s^{j+1} \xi\|_{L^2(R^4 \dR)}^2 + \|\partial_s^j \partial_R \xi\|_{L^2(R^4 \dR)}^2 + \|\partial_s^j \xi\|_{L^2(R^4 \dR)}^2\big) - C \cK_j(s)^2,
\qquad \ell = 1,2.
\end{equation}
The quantity $\cK_j(s)$ will consist of terms that are at least quadratic in $\xi$ and its derivatives, see \eqref{e:Kj} below.  

Substituting \eqref{e:Fjellupper} and \eqref{e:Fjellcoercive} into \eqref{e:xisjIP} and summing over $j$ yields, for all $0\le \widetilde{s}<s\le T$ and some $C,\kappa_0, \kappa_1>0$,
\begin{equation}
\cE_{\temp}^{\le J}(s) +  \kappa_0 \int_{\widetilde{s}}^s \cE_{\temp}^{\le J}(\tau) \dtau \le \kappa_1 \cE_{\temp}^{\le J}(\widetilde{s}) + C \sum_{j=0}^{J-1} \bigg[ \cK_j(s)^2 + \int_{\widetilde{s}}^s \cR_j(\tau)\dtau + \int_{\widetilde{s}}^s \cK_j(\tau)^2\dtau \bigg].
\end{equation}
And then, by Lemma \ref{l:EnEquiv}, we have 
\begin{equation}
\label{e:Etotprimeest}
\cE_{\tot}^{\le J}(s) +  \kappa_0 \int_{\widetilde{s}}^s \cE_{\tot}^{\le J}(\tau) \dtau \le \kappa_1' \cE_{\tot}^{\le J}(\widetilde{s}) + C' \sum_{j=0}^{J-1} \bigg[ \cK_j(s)^2 + \int_{\widetilde{s}}^s \cR_j(\tau)\dtau + \int_{\widetilde{s}}^s \cK_j(\tau)^2\dtau \bigg].
\end{equation}
With \eqref{e:Etotprimeest} in hand, it then remains to estimate the nonlinear terms in the sum on the right side of this inequality; this is the second main technical step mentioned above.  

The remainder of this section is dedicated to the coercivity estimates \eqref{e:Fjellcoercive} and the estimates of the nonlinear terms in~\eqref{e:Etotprimeest}.
The coercivity estimates are treated in  Section~\ref{ss:Coercivity}, while we carry out the estimates of the nonlinear terms in Section~\ref{ss:Nonlinearest}.  Note that we occasionally use the spaces $\cH^k(R^4\dR)$ ($k=1,2$) defined in Remark \ref{r:R4}.

\subsection{Coercivity estimates}

\label{ss:Coercivity}

We turn our attention to the claimed coercivity estimates.  We treat only $\cF_{j,1}$ explicitly; the quantity $\cF_{j,2}$ obeys a lower bound of exactly the same form, by a nearly identical proof, which we omit.  
There are two key ingredients.  The first is a lower bound on the spectral gap of the operator $\cL_{\xi}$, which will allow us to control the mean-zero part of $\xi^{(j)}$.  The Lemma below provides the needed statement.  The second ingredient is the energy identity \eqref{e:Energys}, which allows us to control the mean of $\xi^{(j)}$; we will elaborate on its role below.

\begin{LEMMA}
\label{l:coercivity}
Assume $g\in C^{2}((0,\infty))$, with $g(1) = 1$, $g'(1) = -\frac13$, and $g''(1)>0$. Fix $\e_3\in (0,1)$ such that $g''(y)\ge \frac12 g''(1)$ whenever $|y-1|\le \e_3$.  There exists a constant $\nu>0$, depending only on $g''(1)$ and $\brho$, such that the following statement holds.

Let $\varphi$ be a classical solution of \eqref{e:separated} (with associated parameters $\mu$ and $\brho$) such that  \mbox{$|\varphi''(R)|\le R\e_3$},  $R\in (0,1)$; let $\xi\in C^1([0,1])$ satisfy 
$\|\xi\|_{L^\infty}\le \frac12$ and  $\|\xi_R\|_{L^\infty} \le \frac14\e_3$. Then for any $u\in \cH^1(R^4\dR)$ satisfying $\langle u, 1 \rangle_\xi = 0$, we have 
 \begin{equation} 
 \label{e:Lmucoercive}
    \langle \cL_{\xi} u, u\rangle_\xi \ge \nu (\|u\|_\xi^2 + \|u_R\|_\xi^2).
\end{equation}
The constant $\nu$ may be taken as $\nu = \widetilde{c}\brho^{\frac13} g''(1)$, where $\widetilde{c}$ is an absolute constant.  
\end{LEMMA}
\begin{proof}
We note first of all that, under the assumptions above on $\f$ and $\xi$, we have 
\[
\d_z g'(y_*) = \int_0^1 g''(y_* + \s z) \,\dd\s  \ge \frac12 g''(1).
\]
Indeed, Theorem \ref{t:background} guarantees that $|y_* - 1|\le \frac12 \e_3$, so that $|(y_* + \sigma z)-1| \le |y_* - 1| + 2|\xi_R| \le \e_3$ for all $\s\in [0,1]$.  Therefore, using $\lambda\le \frac43$ and $1+\xi\le \frac32$, we may write
\[
\langle \cL_{\xi} u, u\rangle_\xi
= \brho^{\frac13} \int_0^1 \frac{\delta_z g'(y_*)}{\lambda (1 + \xi)^2} |u_R|^2 R^4 \dR \ge \frac16\brho^{\frac13} g''(1)\int_0^1 |u_R|^2 R^4 \dR.
\]
The classical Poincar\'e inequality, applied to radial functions on the unit ball in $\mathbb{R}^5$, implies that the quantity on the right is bounded below by a constant multiple of $\|u-\overline{u}\|_{L^2(R^4\dR)}^2$, where $\overline{u} = 5 \int_0^1 u(R) R^4 \dR$ is the average of $u$ with the weight $R^4$.  But then 
\[
\|u\|_\xi^2 
\le \|u - \overline{u}\|_\xi^2 
\lesssim 
\|u - \overline{u}\|_{L^2(R^4\dR)}^2
\lesssim \|u_R\|_{L^2(R^4\dR)}^2, 
\]
both implied constants being absolute here.  Since the $L^2(R^4 \dR)$ norm is equivalent to $\|\cdot\|_\xi$, the desired statement follows immediately.
\end{proof}

With the above lemma in hand, we may now specify the parameters $\e_0>0$ and $c_0>0$ in the statement of Theorem \ref{t:Etotest}.  In particular, we choose $\e_0 = \min\{\e_1, \e_2, \e_3\}$, where $\e_1$, $\e_2$, and $\e_3$ are as in Theorem~\ref{t:background}, Lemma~\ref{l:EnEquiv}, and Lemma~\ref{l:coercivity}, respectively.  In what follows, we assume without loss of generality that $\theta_0\le \min\{\theta_2, \eta\e_0^2\}$, where $\theta_2$ is as in Lemma \ref{l:EnEquiv} and $\eta>0$ is an absolute constant small enough to guarantee that $\|\xi\|_{L^\infty}\le \frac12$ and $\|\xi_R\|_{L^\infty} \le \frac14 \e_3$, so that the conclusions of both Lemmas \ref{l:EnEquiv} and \ref{l:coercivity} hold.  
To this end, with $\nu = \widetilde{c} \brho^{\frac13} g''(1)$ as in the statement of Lemma \ref{l:coercivity}, we may take $c_0 = \frac{\widetilde{c}}{6}$, so that the allowed range of $\mu$ is $(-c_0\brho^{\frac13} g''(1), 0) = (-\frac16 \nu, 0)$.

We now project  $\xi^{(j)}$ onto the space of (weighted) mean-zero functions, writing 
\begin{equation} 
\label{e:projection}
\xi^{(j)}=\widetilde{\xi^j} + \overline{\xi^j},
\qquad\quad  
\overline{\xi^j} = \frac{\langle \xi^{(j)}, 1\rangle_\xi}{\|1\|_\xi^2}. 
\end{equation} 

Applying \eqref{e:Lmucoercive} with the choices specified above then yields  
\begin{equation} 
\label{e:Lmuxi3mulwrbd2}
\begin{split} 
\langle (\cL_{\xi} + 3\mu) \xi^{(j)}, \xi^{(j)}\rangle_\xi
& = \langle \cL_{\xi} (\widetilde{\xi^j}), \widetilde{\xi^j}\rangle_\xi
- 3|\mu| \|\widetilde{\xi^j}\|_\xi^2 - 3|\mu| (\overline{\xi^j})^2 \|1\|_\xi^2 \\
& \ge \nu \bigg(\frac12 \|\widetilde{\xi^j}\|_\xi^2  + \|\xi^{(j)}_R\|_\xi^2\bigg)  - 3 |\mu| (\overline{\xi^j})^2 \|1\|_\xi^2.
\end{split}
\end{equation} 
The last term on the right side of this computation makes clear that we will need to recruit a  positive multiple of $(\overline{\xi^j})^2$ from elsewhere in our estimate of $\cF_{j,1}$, in order to achieve a lower bound involving full $L^2$ norm of $\xi^{(j)}$.  Equation \eqref{e:Energys} provides the mechanism, by relating $\overline{\xi^j}$ and $\overline{\xi^{j+1}}$.  More specifically,  differentiating \eqref{e:Energys} with respect to $s$ and rewriting the inner products in terms of $\langle \cdot, \cdot \rangle_\xi$, we  obtain  
\begin{equation} 
\label{e:Kj}
\frac32 b \langle \xi^{(j)}, 1\rangle_\xi -  \langle \xi^{(j+1)}, 1 \rangle_\xi = \cK_j(s):= \frac{1}{b}\partial_s^j \cJ(s) + \int_0^1 \bigg( \frac32 b \xi^{(j)} - \xi^{(j+1)} \bigg) R^4 \lambda^2 \xi \dR.
\end{equation} 
We may therefore write 
\[
(\overline{\xi^{j+1}})^2 \|1\|_\xi^4 = 2b^2 (\overline{\xi^{j}})^2\| 1\|_\xi^4 - 8\cK_j(s)^2 + \bigg( \frac12 b \overline{\xi^j} \|1\|_\xi^2 - 3 \cK_j(s)\bigg)^2
\]
and estimate as follows.
\begin{align*} 
\cF_{j,1}(s) 
& = \|\xi^{(j+1)}\|_\xi^2 + \frac{b^2 k_0}{4} \|\xi^{(j)}\|_\xi^2 - b k_0 \langle \xi^{(j+1)}, \xi^{(j)}\rangle_\xi + \langle (\cL_{\xi} + 3\mu)(\xi^{(j)}), \xi^{(j)} \rangle_\xi \\
& \ge ( 1 - k_0 ) \|\xi^{(j+1)}\|_\xi^2  + \langle (\cL_{\xi} + 3\mu) \xi^{(j)}, \xi^{(j)}\rangle_\xi \\
& \ge (1-k_0-k_1) \|\xi^{(j+1)}\|_\xi^2 + \bigg[ k_1 (\overline{\xi^{j+1}})^2  -\frac32 b^2  (\overline{\xi^j})^2 \bigg]  \|1\|_\xi^2 
+  \nu \bigg(\frac12\|\widetilde{\xi^j}\|_\xi^2 
+ \|\xi^{(j)}_R\|_\xi^2\bigg)\\
& \ge (1-k_0-k_1) \|\xi^{(j+1)}\|_\xi^2 + \min\bigg\{ 2b^2 \Big( k_1 -\frac34 \Big), \frac{\nu}{2} \bigg\} \|\xi^{(j)}\|_\xi^2
+ \nu \|\xi^{(j)}_R\|_\xi^2 - \frac{8k_1}{\|1\|_\xi^2} \cK_j(s)^2.
\end{align*} 
Picking $k_1\in (\frac34, 1)$ and then $k_0\in (0, 1-k_1)$, we obtain \eqref{e:Fjellcoercive}.

\subsection{Estimates of the nonlinear terms}

\label{ss:Nonlinearest}

The final step of our proof of Theorem \ref{t:Etotest} is to establish control over the quantities in the sum on the right side of \eqref{e:Etotprimeest}, namely $\cK_j$ and $\cR_j$ (as defined in \eqref{e:Kj} and \eqref{e:Rj}, respectively), for $j\in \{0, 1, \ldots, J-1\}$.  To clarify the division of labor below, we split the quantity $\cR_j$ as $\cR_j = \cR_j^1 + \cR_j^2 - \cR_j^3$, where 
\begin{small} 
\[
\cR_j^1 =  
\brho^{\frac13}\int_0^1 \partial_s \bigg( \frac{\delta_z g'(y_*)}{\lambda(1+\xi)^2} \bigg) |\partial_s^j \partial_R \xi|^2 R^4 \dR + \int_0^1 \bigg [ |\partial_s^{j+1}\xi|^2 - \frac{b^2}{4} ( 6 - k_0 ) |\partial_s^j \xi|^2 - b k_0 \partial_s^{j+1}\xi\, \partial_s^j \xi  \bigg] \lambda^2 R^4 \xi_s \dR,
\]
\end{small} 
\[
\cR_j^2 = \big\langle \partial_s^j \cN(\xi), 2\partial_s^{j+1}\xi - b k_0 \partial_s^j \xi \big\rangle_\xi,
\qquad 
\cR_j^3 = \big\langle [\partial_s^j, \cL_{\xi}](\xi), 2\partial_s^{j+1}\xi - b k_0 \partial_s^j \xi \big\rangle_\xi,
\]
and we deal with each of these quantities in turn. We make frequent use of the calculus developed in Section \ref{ss:Preliminaries}, especially the principle from Proposition \ref{p:zderivs} that $z = \frac{R\xi_R}{1+\xi}$ should be treated at the same level as $\xi$ for purposes of $L^\infty$ estimates (though not for purposes of $L^2$ estimates).  With these ideas in hand, most of the terms involved in the expansions of $\cK_j$ and $\cR_j$ can be estimated using simple H\"older splits, with $L^2(R^4\dR)$ norms on the terms where the highest-order derivatives appear and $L^\infty$ elsewhere.  However, there are two exceptions.  First, and most importantly, the commutator in $\cR_{J-1}^3$ generates one term where $J+1$ derivatives fall on $\xi$.  We need to integrate by parts to relieve this overloaded term, while still controlling the error terms that this process generates.  And for $J=5$, there are a few instances where a naive $L^2$-$L^2$-$L^\infty$ split is insufficient, and we leverage the Sobolev embeddings $\cH^1(R^4 \dR)\hookrightarrow L^{\frac{10}{3}}(R^4\dR)$ and $\cH^2(R^4 \dR)\hookrightarrow L^{10}(R^4\dR)$ (mentioned earlier in Remark \ref{r:R4}). 

\subsubsection{Estimates on $\cK_j$}
The integral in \eqref{e:Kj} is easy to control; indeed, for $0\le j\le J-1$, we have 
\[
\bigg| \int_0^1 \bigg( \frac32 b \xi^{(j)} - \xi^{(j+1)} \bigg) R^4 \lambda^2 \xi \dR\bigg| \lesssim \|\xi\|_{L^2(R^4 \dR)} \big( \|\partial_s^j \xi\|_{L^2(R^4 \dR)} + \|\partial_s^{j+1} \xi\|_{L^2(R^4 \dR)}\big) \lesssim \cE_{\tot}^{\le J}.
\]
Next, we recall the definition \eqref{e:Js} of $\cJ(s)$:
\[
\cJ(s) = - \int_0^1 \bigg( \frac12 |\xi_s - b\xi|^2 + \frac{\mu \xi^2}{1 + \xi} + (1 + \xi)\xi_R^2\,\cN_2 \bigg) R^4 \l^2 \dR.
\]
Here $\cN_2$, defined in \eqref{e:N2}, is a $J$ times continuously differentiable function of $\lambda$, $y_*$, $\xi$ and $z = \frac{R\xi_R}{1+\xi}$.  Notice that each term inside the integral is quadratic in some combination of $\xi$, $\xi_s$, and $\xi_R$.  Therefore, when $\cJ(s)$ is differentiated $j\in \{0, \ldots, J-1\}$ times with respect to $s$, the result is a sum of terms which are at least quadratic in $\xi$ and its derivatives, multiplied by a bounded function and the weight $R^4$.  The maximum number of derivatives that can appear on a single $\xi$ is $j+1\le J$; the second highest number is $\le \lfloor \frac{j}{2} + 1\rfloor\le J$; and the third highest number (if it appears) is $\lfloor \frac{j}{3} + 1\rfloor$, which is at most $J-3$ as long as $J\ge 5$.  We therefore put $L^2(R^4\dR)$ norms on the highest and second-highest order derivatives, and $L^\infty$ norms on any other quantities involving $\xi$ or $z$.  All in all, we have that 
\[
|\partial_s^j \cJ|\lesssim \cE_{\tot}^{\le J}.
\]
Putting together the last two estimates yields 
\[
\sum_{j=0}^{J-1} \bigg[ \cK_j(s)^2 + \int_{\widetilde{s}}^s \cK_j(\tau)^2\dtau \bigg]
\lesssim (\cE_{\tot}^{\le J}(s))^2 + \int_{\widetilde{s}}^s (\cE_{\tot}^{\le J}(\tau))^2\dtau. 
\]

\subsubsection{Estimates on $\cR_j^1$} This term can also be treated easily.  
For $j\in \{0,\ldots, J-1\}$ and $J\ge 5$, we have 
\begin{align*}
\bigg| \brho^{\frac13} \int_0^1 \partial_s \bigg( \frac{\delta_z g'(y_*)}{\lambda(1+\xi)^2} \bigg) |\partial_s^j \partial_R \xi|^2 R^4 \dR \bigg| 
& \lesssim (\|\xi_s\|_{L^\infty} + \|z_s\|_{L^\infty}) \|\partial_s^j \partial_R \xi\|_{L^2(R^4\dR)}^2  \lesssim (\cE_{\tot}^{\le J})^{\frac32},
\end{align*}

\begin{align*}
& \bigg| \int_0^1 \bigg [ |\partial_s^{j+1}\xi|^2 - \frac{b^2}{4} ( 6 - k_0 ) |\partial_s^j \xi|^2 - b k_0 \partial_s^{j+1}\xi\, \partial_s^j \xi  \bigg] \lambda^2 R^4 \xi_s \dR\bigg| \\
& \hspace{50 mm} \lesssim  \|\xi_s\|_{L^\infty} \big( \|\partial_s^j \xi\|_{L^2(R^4 \dR)}^2 + \|\partial_s^{j+1} \xi\|_{L^2(R^4 \dR)}^2\big) \lesssim (\cE_{\tot}^{\le J})^{\frac32}.
\end{align*}

\subsubsection{Estimates on $\cR_j^2$}
We turn to the estimate of $\cR_j^2$.  The $\cN_1$ contribution is not an issue; the expression $\partial_s^j (\xi^2 \cN_1)$ is at least quadratic in $\xi$ and its time derivatives, with the top-order derivative of order at most $J-1$ and the second-highest order at most $\lfloor \frac{J-1}{2}\rfloor\le J-3$, provided that $J\ge 5$. Therefore
\[
\big| \langle \partial_s^j (\xi^2 \cN_1), 2\partial_s^{j+1}\xi - bk_0 \partial_s^j \xi \rangle_\xi\big| \lesssim\|\partial_s^j (\xi^2 \cN_1)\|_{L^2(R^4 \dR)}  \|2\partial_s^{j+1}\xi - bk_0 \partial_s^j \xi\|_{L^2(R^4\dR)}\lesssim  (\cE_{\tot}^{\le J})^{\frac32}.
\]  
Most of the terms in the expansion of $\partial_s^j(\xi_R^2 \cN_2)$ can be treated similarly to those of $\partial_s^j(\xi^2 \cN_1)$, but the second-highest order derivative can now be as large as $\lfloor \frac{J+1}{2} \rfloor$.  This is still $\le J-3$ for $J\ge 6$, but for $J=5$, an $L^2$-$L^2$-$L^\infty$ split is too crude for one of the terms in the expansion.  Therefore, we leverage the Sobolev embedding $\cH^2(R^4\dR)\hookrightarrow L^{10}(R^4\dR)$, as we previewed earlier, and we write 
\[
\int_0^1 |\xi_{Rss}|^2 |\xi_{sssss}| R^4 \dR
\le \|\xi_{Rss}\|_{\cH^2(R^4 \dR)}^2 \|\xi_{sssss}\|_{L^2(R^4\dR)} \lesssim (\cE_\tot^{\le 5})^{\frac32}. 
\]
With this term out of the way, we can estimate all other contributions easily to obtain 
\[
\big| \langle \partial_s^j (\xi_R^2 \cN_2), 2\partial_s^{j+1}\xi - bk_0 \partial_s^j \xi \rangle_\xi\big| \lesssim (\cE_{\tot}^{\le J})^{\frac32}.
\]
\subsubsection{Estimates on $\cR_j^3$}
We turn now to the commutator $[\partial_s^j, \cL_{\xi}](\xi)$, which requires the most care.  We start by expanding it.  For $j=0$, it is zero; for $j\ge 1$, it takes the form
\begin{equation}
[\partial_s^j, \cL_{\xi}](\xi) 
= \sum_{\substack{j_1 + j_2 + j_3 = j \\ 0\le j_1, j_2; \; 0\le j_3 <j}} a_{j_1 j_2 j_3} (\partial_s^{j_1}\alpha)  D_R\big( \partial_s^{j_2} \beta\; \partial_s^{j_3} \xi_R \big),
\qquad j\ge 1,
\end{equation}
where 
\[
\alpha = - \frac{\brho^{\frac13}}{\lambda^2 (1+\xi)},
\qquad 
\beta = \frac{\delta_z g'(y_*)}{\lambda(1+\xi)^2}.
\]
We estimate the pairing of each term in the sum against $2\partial_s^{j+1} \xi - bk_0 \partial_s^j \xi$.  For most of the terms, we can simply estimate the two parts of the pairing separately via Cauchy--Schwarz, but we will also treat one problematic term below where an integration by parts is necessary.  

We start with the easy cases. First of all, note that $\|2\partial_s^{j+1} \xi - bk_0 \partial_s^j \xi\|_\xi\lesssim (\cE_\tot^{\le J})^{\frac12}$.  
Now, suppose that $j_1 = J-1$.  Since $1\le j\le J-1$, this forces $j_2 = j_3 = 0$.  Then for $J\ge 5$, it is clear that 
\[
\|(\partial_s^{J-1}\alpha)  D_R\big( \beta \xi_R \big)\|_{\xi} \lesssim \|\partial_s^{J-1}\alpha\|_{L^2(R^4\dR)} \|D_R(\beta \xi_R)\|_{L^\infty} \lesssim \cE_{\tot}^{\le J}.
\]
If $j_1 = J-2$, then $(j_2, j_3)$ may be equal to $(0,0)$, $(1,0)$, or $(0,1)$.  The $(j_2, j_3) = (0,0)$ case is dealt with exactly as in the estimate above.  If $(j_2, j_3) = (1,0)$, then 
\begin{align*}
\|(\partial_s^{J-2}\alpha)  D_R\big( (\partial_s\beta) \xi_R \big)\|_{\xi}
& \lesssim \|\partial_s^{J-2} \alpha \|_{L^2(R^4 \dR)} \big( \|\partial_s \partial_R \beta\|_{L^\infty} \|\xi_R\|_{L^\infty} + \|\partial_s \beta\|_{L^\infty} \|\cD^2 \xi\|_{L^\infty} \big) 
\lesssim (\cE_{\tot}^{\le J})^{\frac32}. 
\end{align*}
On the other hand, if $(j_2, j_3) = (0,1)$, then we write 
\begin{align*}
\|(\partial_s^{J-2}\alpha)  D_R\big( \beta \partial_s \xi_R \big)\|_{\xi}
& \lesssim \|\partial_s^{J-2} \alpha \|_{L^2(R^4 \dR)} \|\partial_R \beta\|_{L^\infty} \|\xi_{sR}\|_{L^\infty} \\
& \quad + \|\partial_s^{J-2} \alpha \|_{\cH^2(R^4 \dR)} \|\beta\|_{L^\infty} \|\cD^2 \xi_s\|_{\cH^2(R^4\dR)}  
\lesssim \cE_{\tot}^{\le J}.
\end{align*}
If $1\le j_1\le J-3$, or if $j_1 = 0$ and $1\le j \le J-2$, then we can always put the $L^\infty$ norm on $\partial_s^{j_1}\alpha$, and we only need to control $\|\partial_s^{j_2} \partial_R \beta \; \partial_s^{j_3} \xi_R\|_{L^2(R^4\dR)}$ and $\|\partial_s^{j_2} \beta \; \partial_s^{j_3} \cD^2 \xi\|_{L^2(R^4\dR)}$.  The first of these is easy to deal with, since in this situation, $\max\{j_2, j_3\} \le J-2$ and $\min\{j_2, j_3\}+1\le \lfloor \frac{J}{2} \rfloor \le J-3$, for $J\ge 5$.  Put $L^2(R^4 \dR)$ on the higher-order derivative and $L^\infty$ on the lower one to conclude that 
\[
\|\partial_s^{j_1} \alpha \;\partial_s^{j_2} \partial_R \beta \; \partial_s^{j_3} \xi_R\|_{\xi}\lesssim \cE_{\tot}^{\le J}.
\]
As for the term $\|\partial_s^{j_1} \alpha\;\partial_s^{j_2} \beta \; \partial_s^{j_3} \cD^2 \xi\|_{L^2(R^4\dR)}$, we note that $\max\{j_2, j_3 + 2\} \le J$ (remembering that $j_3 < j$) and $\min\{ j_2, j_3 + 2\} \le \lfloor \frac{J}{2} \rfloor \le J-3$.  So again, we can put the $L^2(R^4 \dR)$ norm on the portion corresponding to $\max\{j_2, j_3 + 2\}$ and put $L^\infty$ on the portion corresponding to $\min\{ j_2, j_3 + 2\}$.  We conclude again that 
\[
\|\partial_s^{j_1} \alpha \;\partial_s^{j_2} \beta \; \partial_s^{j_3} \cD^2 \xi\|_{\xi} \lesssim \cE_{\tot}^{\le J}.
\]
Now we consider the case where $j_1 = 0$ and $j = J-1$.  If $j_2< J-1$, then we can estimate as above for both $\alpha \; \partial_s^{j_2} \partial_R \beta \;\partial_s^{j_3} \partial_R \xi$ and $\alpha \; \partial_s^{j_2} \beta \;\partial_s^{j_3} \cD^2 \xi$: If $J\ge 6$, then we put the $L^2(R^4\dR)$ norm on whichever of $\beta$, $\xi$ carries more derivatives (as determined by the values of $j_2$ and $j_3$), and put $L^\infty$ on everything else.  If $J = 5$, this still works except when a total of three derivatives falls on both $\beta$ and $\xi$, which occurs when $(j_2, j_3)$ is either $(2,2)$ or $(3,1)$.  The corresponding terms are estimated using the Sobolev embeddings mentioned earlier in this section:
\[
\|\partial_{ss} \partial_R \beta \;\partial_{ss} \partial_R \xi \|_\xi \lesssim  \|\partial_{ss} \partial_R \beta\|_{\cH^1(R^4\dR)} \|\partial_{ss} \partial_R \xi \|_{\cH^2(R^4\dR)}
\]
\[
\| \partial_s^{3} \beta \;\partial_s \cD^2 \xi\|_\xi \lesssim \|\partial_{sss}\beta \|_{\cH^1(R^4 \dR)} \|\partial_s \cD^2 \xi \|_{\cH^2(R^4\dR)}.
\]
It is clear that $\|\partial_{ss} \partial_R \xi \|_{\cH^2(R^4\dR)}$ and $\|\partial_s \cD^2 \xi \|_{\cH^2(R^4\dR)}$ are controlled by $(\cE_{\tot}^{\le 5})^{\frac12}$ (cf. Remark \ref{r:R4}).  The contributions involving $\beta$ worth discussing explicitly are the $L^2(R^4\dR)$ norms of $\partial_{ss} \partial_{RR} \beta$ and $\partial_{sss}\partial_R \beta$. Expanding either expression will distribute at most four derivatives between $z$ and $\xi$, so there is no overload of derivatives on any individual $\xi$.  We just need to check, therefore, that a further $L^2$-$L^\infty$ split is still feasible. It turns out that this is not an issue: whenever at least two of the derivatives land on $z$, then we can put the $L^2(R^4\dR)$ norm on this part and $L^\infty$ on everything else, since at most two derivatives can then land on $\xi$.  On the other hand, if at most one derivative lands on $z$, then we can put an $L^\infty$ norm on this part, then an $L^2(R^4\dR)$ norm on the factor with the highest-order derivative, and $L^\infty$ on everything else.  Either way, we may conclude that both $\|\partial_{ss} \partial_{R} \beta\|_{\cH^1(R^4\dR)}$ and $\|\partial_{sss}\beta\|_{\cH^1(R^4\dR)}$ are controlled by $(\cE_{\tot}^{\le 5})^{\frac12}$.

We have now treated every case except the one where $j_1 = j_3 = 0$ and $j_2 = J-1$.  Clearly $\|\alpha \partial_s^{J-1} \beta \;\cD^2 \xi\|_\xi \lesssim \cE_{\tot}^{\le J}$, but the term $\alpha \partial_s^{J-1} \partial_R \beta \; \xi_R$ cannot be estimated directly; after all, $\beta$ depends on~$z$, so applying a total of $J$ derivatives to it results in a total of $J+1$ derivatives on $\xi$.  We therefore return to the full inner product $\langle \alpha \partial_s^{J-1} \partial_R \beta \; \xi_R, \; 2\partial_s^J \xi - bk_0\partial_s^{J-1} \xi\rangle_\xi$, which we rewrite as 
\[
-\brho^{\frac13} \int_0^1 \partial_s^{J-1}\partial_R \bigg( \frac{\delta_z g'(y_*)}{\lambda(1+\xi)^2} \bigg) \xi_R \big( 2\partial_s^J \xi - bk_0 \partial_s^{J-1}\xi \big) \,R^4 \dR.
\]
Expanding the derivative $\partial_s^{J-1} \partial_R$ results in many terms, but the only one that results in an overload of derivatives on $\xi$ can be written 
\[
-\brho^{\frac13}\int_0^1 \partial_z \bigg( \frac{\delta_z g'(y_*)}{\lambda(1+\xi)^2} \bigg) \cdot \frac{R \partial_s^{J-1}\xi_{RR}}{1+\xi} \cdot \xi_R \Big( 2\partial_s^J \xi - bk_0 \partial_s^{J-1} \xi\Big) \,R^4 \dR.
\]
(The other terms can be controlled by $(\cE_{\tot}^{\le J})^{\frac32}$.)
We integrate by parts to remove one of the spatial derivatives from $\partial_s^{J-1} \xi_{RR}$.  Again, this results in many terms, but the only one that is not obviously controlled is
\[
\brho^{\frac13}\int_0^1 \partial_z \bigg( \frac{\delta_z g'(y_*)}{\lambda(1+\xi)^2} \bigg) \cdot \frac{R \partial_s^{J-1}\xi_{R}}{1+\xi} \cdot \xi_R \Big( 2\partial_s^J \xi_R - bk_0\partial_s^{J-1} \xi_R\Big) \,R^4 \dR.
\]
At this point, we can see that the portion involving $\partial_s^{J-1} \xi_R$ is unproblematic, as there is no derivative overload there.  We dismiss this portion and focus on $\partial_s^J \xi_R$, which still contains $J+1$ derivatives.  However, the integrand also contains $\partial_s^{J-1} \xi_R$, so together these two terms form a perfect time derivative of an expression that is not overloaded:
\begin{align*}
& \int_0^1 \partial_z \bigg( \frac{\delta_z g'(y_*)}{\lambda(1+\xi)^2} \bigg) \cdot \frac{R \partial_s^{J-1}\xi_{R}}{1+\xi} \cdot \xi_R \cdot 2\partial_s^J \xi_R  \,R^4 \dR \\
& \qquad = \frac{\dd}{\ds} \bigg[ \int_0^1 z \partial_z \bigg( \frac{\delta_z g'(y_*)}{\lambda(1+\xi)^2} \bigg) |\partial_s^{J-1} \xi_R|^2 \,R^4 \dR \bigg] 
-  \int_0^1 \frac{\dd}{\ds} \bigg[ z \partial_z \bigg( \frac{\delta_z g'(y_*)}{\lambda(1+\xi)^2} \bigg)\bigg] |\partial_s^{J-1} \xi_R|^2 \,R^4 \dR. 
\end{align*}
The integral of the quantity on the right, over the time interval $[\widetilde{s},s]$, can be bounded above by a constant times $(\cE_{\tot}^{\le J}(s))^{\frac32} + \int_{\widetilde{s}}^s (\cE_{\tot}^{\le J}(\tau))^{\frac32}\dtau$. 

\subsubsection{Assembling the estimates}
Assembling all the estimates of this section and substituting them into \eqref{e:Etotprimeest} yields 
\begin{equation}
\cE_{\tot}^{\le J}(s) +  \kappa_0 \int_{\widetilde{s}}^s \cE_{\tot}^{\le J}(\tau) \dtau \le \kappa_1' \cE_{\tot}^{\le J}(\widetilde{s}) + C'' \bigg[ (\cE_{\tot}^{\le J}(s))^{\frac32} + \int_{\widetilde{s}}^s (\cE_{\tot}^{\le J}(\tau))^{\frac32}\dtau\bigg].
\end{equation}
Recalling that $\sup_{\tau\in [0,T]} \cE_{\tot}^{\le J}(\tau)\le \theta_0$ by assumption, we choose $\theta_0$ small enough so that the term in brackets can be absorbed into the left side of the equation, yielding (for some constants $\k_0'', \k_1''>0$)
\begin{equation}
\cE_{\tot}^{\le J}(s) +  \kappa_0'' \int_{\widetilde{s}}^s \cE_{\tot}^{\le J}(\tau) \dtau \le \kappa_1''\cE_{\tot}^{\le J}(\widetilde{s}).
\end{equation}

\appendix
\renewcommand{\thesection}{A}
\setcounter{subsection}{0}
\section*{Appendix}

\label{s:appendix}

\textbf{\underline{Proof of Proposition \ref{p:productrule}}.} 
We remind the reader that we wish to prove that for any $A\in \cP_k$, there exist constants $c^{ABC}$ such that 
\begin{equation}
	A(fg) = \sum_{\ell=0}^k \sum_{B\in \cP_\ell, \; C\in \cP_{k - \ell}} c^{ABC} Bf \cdot Cg.
\end{equation}   

If $k$ is $0$ or $1$, then there is nothing to prove.  Assume the statement has been shown for $k = j$, and that $A\in \cP_{j+1}$.  Suppose first that $j$ is even, so $A = \p_R A'$ for some $A'\in \cP_j$, and 
\[
A'(fg) = \sum_{\ell=0}^j \sum_{B'\in \cP_\ell, \; C'\in \cP_{j - \ell}} c^{A'B'C'} B'f \cdot C'g.
\]
We apply $\p_R$ to each term $B' f \cdot C' g$.  Note that $\ell$ and $j-\ell$ are either both even or both odd.  In the first case, we have $\p_R (B'f \cdot C'g) = \p_R B' f \cdot C' g + B' f \cdot \p_R C' g$, which is a linear combination of terms of the desired form.  If $\ell$ and $j-\ell$ are both odd, then we write 
\begin{align*}
	\p_R (B' f \cdot C'g)
	& = \big(D_R - \frac4R\big) B'f \cdot C' g + B'f \big(D_R - \frac4R\big)C'g \\
	& = D_R B' f \cdot C'g - \frac8R B'f \cdot C'g + B' f \cdot D_R C' g.
\end{align*}
This is once again a linear combination of terms of the desired form.

The case where $j$ is odd is even easier.  In this case we take $A' \in \cP_j$, expanded just as before, except that now exactly one of $\ell$, $j-\ell$ is even.  Suppose $\ell$ is even.  Then 
\[
\frac1R (B' f \cdot C'g) = B' f \cdot \frac1R C' g,
\]
\[
D_R (B' f \cdot C'g) = \p_R B' f \cdot C' g + B' f \cdot D_R C' g.
\]
Thus, in all cases, $A(fg)$ has the desired form; this completes the proof.

\medskip
\textbf{\underline{Proof of Proposition \ref{p:Pkcontrol}.}} We wish to establish estimates of the form $\|AX\| \lesssim \|\cD^k X\|$, for $A\in \cP_k$ and $\|\cdot\|$ either the $L^2(R^4\dR)$ norm or the $L^\infty$ norm. If $k$ is $0$ or $1$, then there is nothing to prove.  We start by proving the desired statements in the case where $A = \cD^j(\frac1R \partial_R \cdot)$, $0\le j\le k-2$.  Denote $Y_j = \cD^j\big( \frac1R \partial_R X\big)$, and note that Lemma \ref{l:Dkcommutator} yields the identity 
\[
R^{j+4} \cD^{j+2} X = R^{j+4} \big( R\partial_R + j + 5 \big) Y_j = \partial_R ( R^{j+5} Y_j ), 
\]
which we can integrate to obtain 
\[
R^{\frac52} Y_j
= \frac{1}{R^{j+\frac52}} \int_0^R r^{j+4} \cD^{j+2} X(r)\dr + \frac{C_j}{R^{j+\frac52}}.
\]
Using the Cauchy--Schwarz inequality and the finiteness of $\|\cD^{j+2}X\|_{L^2(R^4\dR)}$ (which holds under the hypotheses of either (a) or (b)), we see that the first term on the right is $o(1)$ as $R\to 0+$.  Therefore, $R^{\frac52} Y_j = o(1)$ if and only if $C_j=0$.  We prove that $C_j=0$ by induction. If $C_0\ne 0$, then we have $|Y_0| = \frac1R |\partial_R X| \gtrsim \frac{|C_0|}{R^5}$ near zero (with $Y_0$ having fixed sign) and therefore $|\partial_R X|\gtrsim \frac{|C_0|}{R^4}$, which is impossible, since $\partial_R X\in L^2(R^4\dR)$ by assumption.  Therefore $C_0 = 0$.  Now assume that $C_{j-1}$ is zero for some $j\in \{1, \ldots, k-2\}$.  If $C_j \ne 0$, then $|Y_j| \gtrsim \frac{|C_j|}{R^{j+5}}$ near zero (with $Y_j$ having fixed sign).  The precise consequence for $Y_{j-1}$ depends on whether $j$ is even or odd, but in either case, one can show that $|Y_{j-1}| \gtrsim \frac{|C_j|}{R^{j+4}}$, which means that $R^{\frac52} Y_{j-1}$ cannot be $o(1)$ and $C_{j-1}$ cannot be zero, a contradiction.  

The argument above allows us to write 
\begin{equation} 
\label{e:Yjformula}
\cD^j\bigg(\frac1R\partial_R X\bigg) = \frac{1}{R^{j+5}} \int_0^R r^{j+4} \cD^{j+2} X(r)\dr, 
\end{equation} 
so that 
\begin{align*}
\int_0^1 \bigg|\cD^j\bigg(\frac1R\partial_R X\bigg)\bigg|^2 R^4\,\dR 
& \le \int_0^1 \bigg( \frac{1}{R^{2j+5}}  \int_0^R r^{2j+8} |\cD^{j+2} X(r)|^2 \dr\bigg) \dR \\
& = \int_0^1 \bigg( \int_r^1 \frac{\dR}{R^{2j+5}} \bigg) r^{2j+4} \cdot |\cD^{j+2} X(r)|^2 r^4\dr \lesssim \int_0^1 |\cD^{j+2} X(r)|^2 r^4\dr.
\end{align*}
Even more simply, 
\[
\bigg\| \cD^j \bigg( \frac1R \partial_R X \bigg) \bigg\|_{L^\infty} \le \frac{1}{j+5} \|\cD^{j+2} X\|_{L^\infty}.
\]

We have now proved Proposition \ref{p:Pkcontrol} in its entirety when $k\in \{0,1,2\}$, and also for $k>2$ in the special case where $A$ takes the form $A = \cD^{k-2} \big( \frac1R \p_R \cdot \big)$.  In other words, 
\begin{equation} 
	\label{e:Pkcontrolbootstrap}
    \begin{split} 
	\int_0^1 |\cD^{k-2} \big( \tfrac1R \p_R X\big) |^2 R^4 \dR & \lesssim \int_0^1 |\cD^k X|^2 R^4 \dR, \\
\bigg\| \cD^{k-2} \bigg( \frac1R \partial_R X \bigg) \bigg\|_{L^\infty} & \lesssim \|\cD^{k} X\|_{L^\infty}.
    \end{split} 
	\qquad k\ge 2.
\end{equation} 
In the general case $A\in \cP_k$, $k>2$, we write 
\[
A = \cD^{k_1} \frac1R \p_R \cD^{k_2} \frac1R \p_R \cD^{k_3} \cdots \frac1R \p_R \cD^{k_n}, 
\]
where $k_1 + k_2 + \cdots + k_n + 2(n-1) = k$, and $k_j\ge 0$ for each $j\in \{1, \ldots, n\}$.  (Each of $k_2, \ldots, k_n$ is even, while $k_1$ may be either even or odd.)  We also denote $X_0 = \cD^{k_n} X$ and for $\ell \in \{1, \ldots, n-1\}$, we define 
\[
X_\ell = \cD^{k_{n-\ell}} \frac1R \p_R X_{\ell-1}.
\]
Thus 
\[
AX = X_{n-1} = \cD^{k_1} \frac1R \p_R X_{n-2} = \cD^{k_1} \frac1R \p_R \cD^{k_2} \frac1R \p_R X_{n-3} = \cdots = \cD^{k_1} \frac1R \p_R \cD^{k_2} \frac1R \p_R \cD^{k_3} \cdots \frac1R \p_R X_0.
\]
Then applying \eqref{e:Pkcontrolbootstrap} repeatedly, we get
\begin{align*}
	\int_0^1 |AX|^2 R^4 \dR 
	& = \int_0^1 |\cD^{k_1} \tfrac1R \p_R X_{n-2}|^2 R^4  \dR \\
	& \lesssim \int_0^1 |\cD^{k_1 + 2} X_{n-2}|^2 R^4 \dR = \int_0^1 |\cD^{k_1 + k_2 + 2} \frac1R \p_R X_{n-3}|^2 R^4 \dR \\
	& \lesssim \int_0^1 |\cD^{k_1 + k_2 + 4} X_{n-3}|^2 R^4 \dR \\
	& \lesssim \cdots \lesssim \int_0^1 |\cD^{k_1 + \cdots + k_{n-1} + 2(n-1)} X_0|^2 R^4  \dR \\
	& = \int_0^1 |\cD^k X|^2 R^4 \dR.
\end{align*}
Similarly, 
\[
\|AX\|_{L^\infty} \lesssim \|\cD^{k_1} \tfrac1R \partial_R X_{n-2} \|_{L^\infty} \lesssim \cdots \lesssim \|\cD^{k_1 + \cdots + k_{n-1} + 2(n-1)} X_0\|_{L^\infty} = \|\cD^k X \|_{L^\infty}.
\]
This completes the proof.

\medskip
\textbf{\underline{Proof of Proposition \ref{p:zderivs}}.}
Our goal is to establish the $L^2$ and $L^\infty$ estimates \eqref{e:Azest1} and \eqref{e:Azest2} on $\partial_s^p Az$, where $A\in \cP_k$.  These estimates are trivial when $p=k=0$; we write out the easy case where $p=0$, $k=1$, and $A = \partial_R$ before moving on to the more general argument.  We have 
\begin{equation}
	\label{e:zR}
	z_R = \frac{R D_R \xi_R}{1+\xi} - \frac{3\xi_R}{1+\xi} - \frac{R\xi_R^2}{(1+\xi)^2}.
\end{equation}
Using the Sobolev and Hardy inequalities \eqref{e:Sob} and \eqref{e:Hardy}, we obtain 
\begin{align*}
	\int_0^1 |z_R|^2 R^4 \dR 
	& \lesssim \int_0^1 |\cD^2 \xi|^2 R^4 \dR + \int_0^1 |\xi_R|^2 R^4 \dR + \| R^2 \xi_R \|_{L^\infty}^2 \int_0^1 \bigg| \frac{\xi_R}{R} \bigg|^2 R^4 \dR \\
	& \lesssim \bigg(\sum_{j=1}^2 \int_0^1 |\cD^j \xi|^2 R^4 \dR\bigg)  \cdot \bigg( 1 + \sum_{j=1}^2 \int_0^1 |\cD^j \xi|^2 R^4 \dR \bigg). 
\end{align*}
Furthermore, 
\[
\|z_R\|_{L^\infty}^2 \lesssim \bigg( \sum_{j=1}^4 \int_0^1 |\cD^j \xi|^2 R^4 \dR\bigg)  \cdot \bigg( 1 + \sum_{j=1}^3 \int_0^1 |\cD^j \xi|^2 R^4 \dR \bigg). 
\]

Now we move on to the general case.  
By an elementary induction argument (the proof of which we omit), it can be shown that whenever $A\in \cP_k$ and $p\ge 0$, we have that  $\partial_s^p Az$ is a linear combination of terms of the form 
\begin{equation}
\label{e:Aztype12terms}
\frac{(\partial_s^{p_1} A_1 \xi)(\partial_s^{p_2} A_2 \xi) \cdots (\partial_s^{p_i} A_i \xi)}{(1+\xi)^{i}} 
\qquad \text{ and } \qquad 
\frac{R(\partial_s^{p_1} A_1 \xi)(\partial_s^{p_2} A_2 \xi) \cdots (\partial_s^{p_i} A_i \xi)}{(1+\xi)^{i}},
\end{equation}
which we refer to as ``Type 1'' and ``Type 2'' terms, respectively.  Here  $A_{\ell} \in \cP_{j_\ell}$ for each $\ell$, each $n_\ell := j_\ell + p_\ell$ is positive, and at least one of the $j_\ell$'s is positive.   For Type~1 terms, the index $i$ lies between $1$ and $k+p$, and the sum of the $n_\ell$'s is $k+p$; for Type 2 terms, $i$ lies between $1$ and $k+p+1$, and the $n_\ell$'s sum to $k+p+1$.

We show the estimates \eqref{e:Azest1} and \eqref{e:Azest2} hold when $\partial_s^p Az$ is replaced by either a `Type 1' term $T_1$ or a `Type 2' term $T_2$.  We consider the Type 1 case first and write $T_1 = \frac{(\p_s^{p_1} A_1 \xi)(\p_s^{p_2} A_2 \xi) \cdots (\p_s^{p_i} A_i \xi)}{(1+\xi)^{i}}$, where $1\le i\le k+p$, $A_\ell\in \cP_{j_\ell}$ for each $\ell = 1, \ldots, i$, and $n_1 + \cdots + n_i = k+p$.  Without loss of generality, we assume that $n_1 \ge n_2 \ge \cdots \ge n_i$. 

If $i=1$, then  $T_1 = \frac{\partial_s^{p} A_1 \xi}{1+\xi}$ for some $A_1\in \cP_k$, so that $\|T_1\|_{L^2(R^4\dR)}^2 \lesssim \cE_{\tot}^{\le k+p}$ and $\|T_1\|_{L^\infty}^2 \lesssim \cE_{\tot}^{\le k+p+3}$.  We therefore assume $i\ge 2$ for the rest of our Type 1 estimates; this also forces $k+p\ge 2$.  We consider the case where $i=k+p=2$ individually.  The only possibilities for $T_1$ here are of the form $\frac{X\xi_R}{(1+\xi)^2}$, where $X$ is either $\partial_s \xi$ or  $\xi_R$.  In either case, 
\[
\begin{split} 
\int_0^1 \bigg| \frac{X \xi_R}{(1+\xi)^2} \bigg|^2 R^4 \dR 
& \lesssim  \int_0^1 \bigg| \frac{\xi_R}{R} \bigg|^2 R^4\dR \cdot\|R X\|_{L^\infty}^2 \lesssim \cE_{\tot}^{\le 2} \cdot \cE_{\tot}^{\le 3}
\end{split} 
\]
and 
\[
\bigg\| \frac{X\xi_R}{(1+\xi)^2} \bigg\|_{L^\infty}^2 \lesssim \|\xi_R\|_{L^\infty}^2 \|X\|_{L^\infty}^2 \lesssim \big( \cE_{\tot}^{\le 4} \big)^2.
\]
The remaining cases can be treated together.  When $i\ge 2$ and $k+p\ge 3$, we can write
\begin{align*} 
	\int_0^1 |T_1|^2 R^4 \dR 
	& \lesssim \bigg( \int_0^1 |\partial_s^{p_1} A_1 \xi|^2 R^4 \dR \bigg) \|\partial_s^{p_2} A_2 \xi\|_{L^\infty}^2 \cdots \|\partial_s^{p_i} A_i\xi\|_{L^\infty}^2  \lesssim \cE_{\tot}^{\le n_1} \big( \cE_{\tot}^{\le n_2 + 3} \big)^{i-1} 
\end{align*}
Since $n_1 \le k+p$ and $n_2 \le \lfloor \frac{k+p}{2} \rfloor\le k+p-2$ for $k+p\ge 3$, the above can be bounded above by the right side of \eqref{e:Azest1}.  If we use an $L^\infty$ norm instead of an $L^2$ norm, then we obtain 
\begin{align*} 
	\|T_1\|_{L^\infty}^2 
	& \lesssim \cE_{\tot}^{\le n_1+3} \big( \cE_{\tot}^{\le n_2 + 3} \big)^{i-1} , 
\end{align*}
which can similarly be bounded above by the right side of \eqref{e:Azest2}.

We turn to the Type 2 terms $T_2 = \frac{R(\partial_s^{p_1} A_1 \xi)(\partial_s^{p_2} A_2 \xi) \cdots (\partial_s^{p_i} A_i \xi)}{(1+\xi)^{i}}$, where $1\le i\le k+p+1$, $A_\ell\in \cP_{j_\ell}$ for each $\ell = 1, \ldots, i$, and $n_1 + \cdots + n_i = k+p+1$. As before, we take $n_1 \ge n_2 \ge \cdots \ge n_i$ without loss of generality.  If $i = 1$, then $T_2 = \frac{R \partial_s^{p_1} A_1 \xi}{1+\xi}$, where $A_1\in \cP_{j_1}$ and $p_1 + j_1 = p+k+1$; in this case, we have $\|T_2\|_{L^2(R^4\dR)}^2 \lesssim \cE_{\tot}^{\le k+p+1}$ and $\|T_2\|_{L^\infty}^2 \lesssim \cE_{\tot}^{\le k+p+3}$.  We therefore assume that $i\ge 2$ for the rest of the proof.  

We still need to consider several subcases; the first one we treat is the one where $k+p=1$.  Then $i=2$, which forces $n_1 = n_2 = 1$, and at least one of $j_1, j_2$ must be equal to $1$.  The only possibilities are of the form $T_2 = \frac{RX\xi_R}{(1+\xi)^2}$, where $X$ is either $\partial_s \xi$ or $\xi_R$.  Either way, we have \[
\|T_2\|_{L^2(R^4\dR)}^2 
\lesssim \int_0^1 \bigg| \frac{\xi_R}{R}\bigg|^2 R^4 \dR \cdot \|R^2 X\|_{L^\infty}^2
\lesssim \cE_{\tot}^{\le 2} \cdot \cE_{\tot}^{\le 2},
\]
\[
\|T_2\|_{L^\infty}^2 
\lesssim \|\xi_R \|_{L^\infty}^2\cdot \|RX\|_{L^\infty}^2
\lesssim \cE_{\tot}^{\le 4} \cdot \cE_{\tot}^{\le 3}.
\]

The case where $k+p = 2$, $i=3$, and $n_1 = n_2 = n_3=1$ is similar.  The only possibilities for $T_2$ are of the form $T_2 = \frac{RXY\xi_R}{(1+\xi)^3}$, where $X,Y\in \{\partial_s \xi, \xi_R\}$.  By essentially the same calculation as the one just above, we have $\|T_2\|_{L^2(R^4\dR)}^2 
\lesssim \cE_{\tot}^{\le 2} \cdot (\cE_{\tot}^{\le 3})^2$ and $\|T_2\|_{L^\infty}^2 
\lesssim \cE_{\tot}^{\le 3} \cdot (\cE_{\tot}^{\le 4})^2$ in this case.

In all the cases we have not yet treated, we have $k+p\ge 2$ and $i\ge 2$.  In this case, we have $n_1\le k+p$ and $n_2 \le \lfloor \frac{k+p+1}{2}\rfloor \le k+p-1$.  If $i\ge 3$, then we must also consider $n_3, \ldots, n_i$; however, we claim that $n_3$ at most $k+p-2$ in all remaining cases (and therefore the same is true of all the other $n_\ell$'s).  When $k+p\ge 3$, this is automatic from the rule $n_3 \le \lfloor \frac{k+p+1}{3}\rfloor$. When $k+p=2$, we have to argue slightly differently: It is possible that $n_3 = 1 = k+p-1$, however, this forces $n_1 = n_2 = 1$ and $i=3$ as well, and we have already explicitly considered this case above.  Therefore, in the computations below, we may take $n_1 \le k+p$, $n_2 \le k+p-1$, and $n_\ell\le k+p-2$ for all $\ell\in \{3, \ldots, i\}$.  In this case, we have (interpreting $\partial_s^{p_\ell} A_\ell \xi$ as $1$ for $\ell > i$, whenever such terms appear below)
\begin{align*} 
\|T_2\|_{L^2(R^4\dR)}^2 
	& \lesssim \bigg( \int_0^1 |\partial_s^{p_1} A_1 \xi|^2 R^4 \dR \bigg) \|R \partial_s^{p_2} A_2 \xi\|_{L^\infty}^2 \cdot \|\partial_s^{p_3} A_3 \xi\|_{L^\infty}^2 \cdots \|\partial_s^{p_i} A_i\xi\|_{L^\infty}^2 \\
	& \lesssim \cE_{\tot}^{\le k+p} \cdot \cE_{\tot}^{\le k+p+1} \cdot 
	\big( \cE_{\tot}^{\le k+p+1} \big)^{i-2} = \cE_{\tot}^{\le k+p} \cdot 
	\big( \cE_{\tot}^{\le k+p+1} \big)^{i-1},
\end{align*}
\begin{align*} 
\|T_2\|_{L^\infty}^2 
	& \lesssim  \|\partial_s^{p_1} A_1 \xi\|_{L^\infty}^2 \| R\partial_s^{p_2} A_2 \xi\|_{L^\infty}^2 \cdot \|\partial_s^{p_3} A_3 \xi\|_{L^\infty}^2 \cdots \|\partial_s^{p_i} A_i\xi\|_{L^\infty}^2 \\
	& \lesssim \cE_{\tot}^{\le k+p+3} \cdot 
	\big( \cE_{\tot}^{\le k+p+1} \big)^{i-1}. 
\end{align*}
This completes the proof.

\section*{Acknowledgments and disclosures}

\textbf{Funding.}  JJ is supported in part by the NSF grants DMS-2306910 and DMS-2605589 and by the Simons Fellowship (Record ID: SFI-MPS-SFM-00021475). TL is supported in part by NSF grant DMS-2408585. 

\medskip
\textbf{AI disclosure statement.}  
During the final stage of the manuscript preparation, 
Claude Opus 5.5 (Anthropic) was used to check the proofs and cross-references throughout the manuscript. The authors have reviewed all contents, including minor corrections and  suggested edits on mathematics and exposition, and have finalized the manuscript themselves. The authors take full responsibility for the manuscript.


\end{document}